\documentclass[smallcondensed,envcountsect]{svjour3}

\usepackage{color}
\newcommand{\modColor}[1]{#1}

\usepackage{geometry}                
\usepackage{graphicx}
\usepackage{amssymb}
\usepackage{amsmath} 
\usepackage{color}

\usepackage{amsthm}
\usepackage{bm}
\usepackage{epstopdf}
\usepackage{accents} 
\usepackage{stmaryrd}
\usepackage{adjustbox}

\definecolor{redcol}{rgb}{1.,0.,0.0} 
\definecolor{lnkcol}{rgb}{0.,0.,0.0} 
\usepackage[pdftex,pdfpagelabels=true,bookmarksdepth=3]{hyperref}
\hypersetup{  colorlinks=true,
              citecolor=lnkcol,%
              linkcolor=lnkcol,%
              urlcolor=lnkcol,
            pdfborder=false,
            bookmarks=true,
            bookmarksnumbered=true,
            pdftitle = {Free-Stream Preservation for Curved Geometrically Non-Conforming Discontinuous Galerkin Spectral Elements},
pdfsubject = {},
pdfauthor = {D.Kopriva,F.Hindenlang,T.Bolemann,G.Gassner},
pdfkeywords = {}
 }
\usepackage[all]{hypcap}  
\usepackage[dvips]{thumbpdf}

\usepackage{cleveref} 

  \crefname{chapter}{Chap.}{Chaps.}       
  \Crefname{chapter}{Chapter}{Chapters}
  \crefname{section}{Sec.}{Secs.}
  \Crefname{section}{Section}{Sections}
  \crefname{table}{Table}{Tables}
  \Crefname{table}{Table}{Tables}
  \crefname{figure}{Fig.}{Figs.}
  \Crefname{figure}{Figure}{Figures}
  \crefname{equation}{}{}
  \Crefname{equation}{Equation}{Equations}
  \crefname{algorithm}{Alg.}{Algs.}
  \Crefname{algorithm}{Algorithm}{Algorithms}
  \crefname{thm}{Thm.}{Thms.}
  \Crefname{thm}{Theorem}{Theorems}
  \crefname{lem}{Lemma}{Lemmas}
  \Crefname{lem}{Lemma}{Lemmas}
  \crefname{cor}{Corr.}{Corrs.}
  \Crefname{cor}{Corrollary}{Corrollaries}
  \crefname{prop}{Prop.}{Props.}
  \Crefname{prop}{Proposition}{Propositions}
  \crefname{rem}{Rem.}{Rems.}
  \Crefname{rem}{Remark}{Remarks}
  \crefname{appendix}{Appendix}{Appendices}
  \Crefname{appendix}{Appendix}{Appendices}
  
\usepackage{stackengine}
\stackMath

\usepackage{accsupp} 

\makeatletter
\def\munderbar#1{\underline{\sbox\tw@{$#1$}\dp\tw@\z@\box\tw@}}
\makeatother

\newcommand\iprod[1]{\left\langle #1\right\rangle} 				

\DeclareMathAccent{\spacevec}{\mathord}{letters}{126}           
\newcommand\acclrvec[1]{\accentset{\,\leftrightarrow}{#1}}	

\newcommand\statevec[1]{\mathbf #1}					
\newcommand\statevecg[1]{\boldsymbol #1}			
\newcommand\contrastatevec[1]{\tilde{\mathbf #1}} 			
\newcommand\bigstatevec[1]{\acclrvec{{\mathbf #1}}}		
\newcommand\bigcontravec[1]{\acclrvec{\tilde{\mathbf #1}}} 	

\newcommand\dS{\,\operatorname{dS} }

\newcommand\ec{\mathrm{ec}}

\newcommand\mmatrix[1]{\underbar{#1}}				
\newcommand{\average}[1]{\left\{\!\!\left\{#1\right\}\!\!\right\}}   

\theoremstyle{plain}
\theoremstyle{remark}
\newtheorem{rem}{Remark}

\usepackage{xargs}                      
\usepackage[pdftex,dvipsnames]{xcolor}  
\usepackage[colorinlistoftodos,prependcaption,textsize=tiny]{todonotes}
\newcommandx{\unsure}[2][1=]{\todo[linecolor=blue,backgroundcolor=blue!25,bordercolor=blue,#1]{#2}}
\newcommandx{\change}[2][1=]{\todo[linecolor=red,backgroundcolor=red!25,bordercolor=red,#1]{#2}}
\newcommandx{\info}[2][1=]{\todo[linecolor=OliveGreen,backgroundcolor=OliveGreen!25,bordercolor=OliveGreen,#1]{#2}}
\newcommandx{\improvement}[2][1=]{\todo[linecolor=Plum,backgroundcolor=Plum!25,bordercolor=Plum,#1]{#2}}
\newcommandx{\thiswillnotshow}[2][1=]{\todo[disable,#1]{#2}}

\title{Free-Stream Preservation for Curved Geometrically Non-Conforming Discontinuous Galerkin Spectral Elements}
\author{David~A.~Kopriva, Florian~Hindenlang, Thomas~Bolemann and Gregor~J.~Gassner}
\date{}                                           
\institute{David A. Kopriva \at  Professor Emeritus, Department of Mathematics, The Florida State University, Tallahassee, FL 32306, USA (\email{kopriva@math.fsu.edu}) \\
           Florian J. Hindenlang \at Max Planck Institute for Plasma Physics, Boltzmannstra{\ss}e 2, D-85748 Garching, Germany \\
Thomas Bolemann \at Institute for Aerodynamics and Gas Dynamics, University of Stuttgart, Pfaffenwaldring 21, Stuttgart, Germany\\ 
Gregor J. Gassner \at Mathematical Institute and Center for Data and Simulation Science, University of Cologne, Cologne,\\ 
           }

\begin{document}
\maketitle
\begin{abstract}
The under integration of the volume terms in the discontinuous Galerkin spectral element approximation introduces errors at non-conforming element faces that do not cancel and lead to free-stream preservation errors. We derive volume and face conditions on the geometry under which a constant state is preserved. From those, we catalog eight constraints on the geometry that preserve a constant state. Numerical examples are presented to illustrate the results.
\end{abstract}
\section{Introduction}

The use of high order methods with naively computed geometry and metric terms can lead to undesired errors. One class of such errors often discussed in the literature are the so-called free-stream preservation errors, where a uniform flow without external forces changes in time due solely to errors in the metric terms. The importance of satisfying free-stream preservation by a numerical approximation method has been emphasized many times e.g. \cite{Thomas&Lombard1979,Visbal2002,Visbal1999,Nonomura2010197,Persson20091585,kopriva2006}. Failure to satisfy the condition can induce artificial sources that can generate spurious waves and affect numerical and flow stability.

In the absence of external forces, inviscid and viscous flows, for instance, are free-stream preserving. On a fixed three-dimensional domain $\Omega\left(\spacevec x\right)$, the Euler equations of gas dynamics and the compressible Navier-Stokes equations can be written in conservation form as
 \begin{equation}
{{{\statevec{u}}}_t} + \nabla  \cdot \bigstatevec f=0,
\label{eq:OrigConsLaw}
\end{equation}
where ${\statevec{u}}$ is the five dimensional state vector and
\begin{equation}\bigstatevec{f}  = \sum\limits_{n = 1}^3 {{{\statevec{f}}_n}{{\hat x}_n}} \end{equation}
 is the (covariant) flux vector, whose components are state vectors.
 We denote the solution and flux vectors here by bold face and spatial vectors by arrows. Double arrows represent space vectors of state vectors. The viscous terms for a constant flow vanish because the gradients are zero. The remaining Euler fluxes are homogeneous functions of degree one, i.e,  $\statevec{f}_{n} = \mmatrix A_{n}\left(\statevec{u}\right)\statevec{u}$. When $\statevec{u}=\statevec{c}=\text{const}$ then $\statevec{f}_{n}=\mmatrix A_{n}\left(\statevec{c}\right)\statevec{c}=\text{const}$. Then the divergence of the constant flux is zero, which leads to the statement that a uniform flow stays uniform for all time.

Equations of the type (\eqref{eq:OrigConsLaw}) can be accurately and efficiently approximated by spectral element methods, especially discontinuous Galerkin versions (DGSEM) \cite{canuto2007},\cite{Hestahven:1008th},\cite{Karniadakis:2005fj},\cite{Kopriva:2009nx}. The physical domain is subdivided into elements and on those elements the solution and fluxes are approximated by high order polynomials. 
The high order of the methods enables them to use curved elements to accurately approximate curved physical boundaries \cite{ISI:A1986A176500008},\cite{Lovgren:2009zr},\cite{ISI:000306588600006}.  If curved elements are used not just at the physical boundaries but also in the neighboring volume, then thinner and longer (anisotropic) elements can be used \modColor{without physical boundaries crossing interior element boundaries} \cite{Hindenlang:2014gl}. The importance of accurately curving boundaries is discussed in \cite{Hindenlang:2014gl} and \cite{NelsonDissertation}. Curved element mesh generation is currently a topic of research \cite{Persson:2009mz},\cite{Johnen:2014eu},\cite{Xie:2013rw},\cite{Moxey:2015xe},\cite{Fortunato2015}.

Hexahedral or quadrilateral versions of spectral element methods have a number of computational advantages  \cite{Kovalev:2005qv}. Most notable is the available tensor product basis, which makes the computation of spatial derivatives very efficient. They are also well suited for computing thin layers like boundary layers.
But hexahedral or quadrilateral spectral element methods also have disadvantages: They are not easily generated and they are not easy to refine locally.

Geometrically, the easiest way to refine a quadrilateral or hexahedral mesh is to make the mesh non-conforming \cite{NME:NME2800}.  A mesh is \emph{conforming} when neighboring elements share a corner point, an entire edge, or an entire face. In geometrically non-conforming meshes, two elements might share a partial edge or face. 
Geometrically non-conforming elements require a modification of the usual method to couple them. A conservative approach was introduced in \cite{Kopriva:1996:JCP96b} and \cite{Koprivaetal2000} where the surface fluxes were computed on a separate mortar space and then projected back onto the element surfaces. 
Stability for straight-sided elements in case of linear hyperbolic systems was shown in \cite{Bui-Tanh:2012zp},\cite{L.-Friedrichs:2017hm}. Entropy stability for non-linear hyperbolic systems is more delicate and needs a different approach \cite{Lucas-Friedrich:2018ng}.
Other examples include static mesh application of the mortar methods \cite{Kopera201492} and sliding mesh applications in \cite{Zhang2015147}.

\begin{figure}[!htbp] 
  \centering
  {\includegraphics[width=0.55\textwidth]{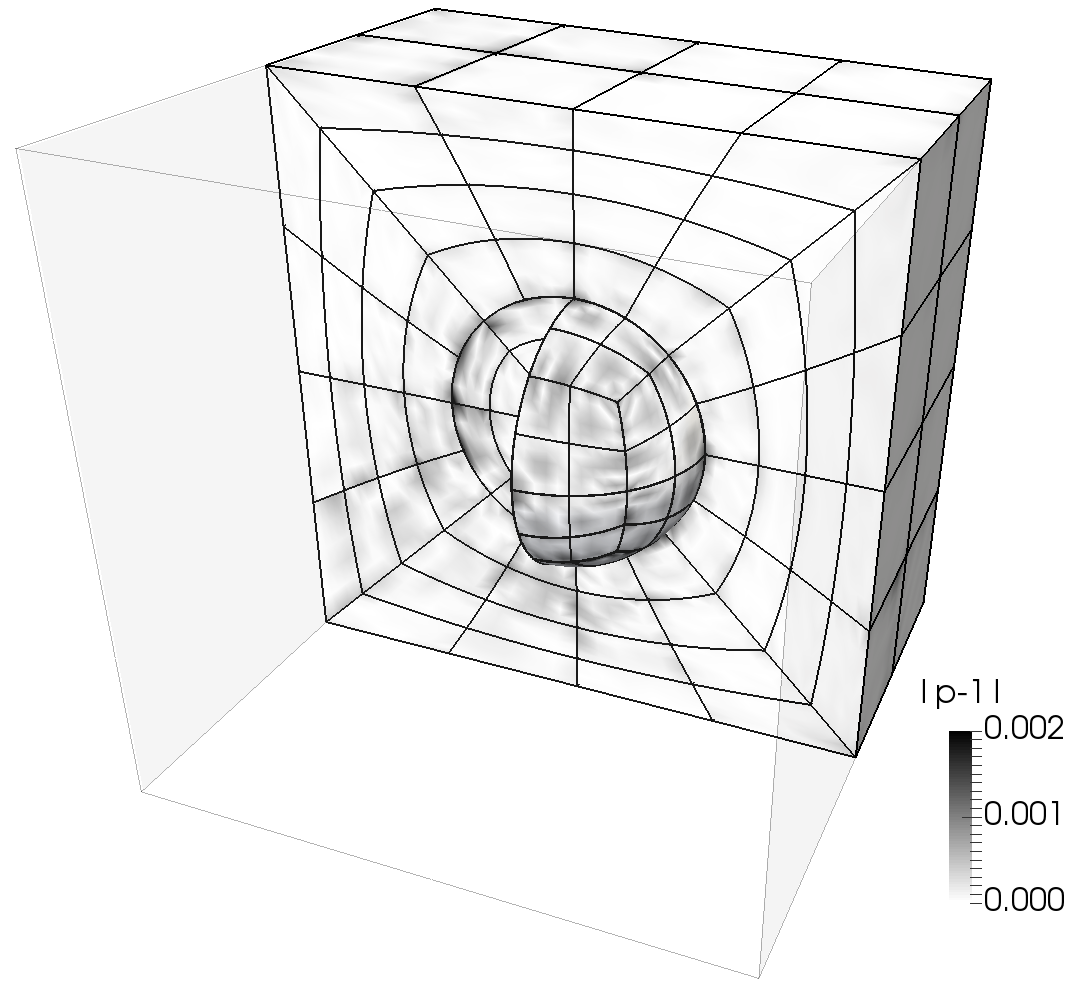} }
   \caption{Cutaway of the absolute value of pressure fluctuations 
   with geometrically induced free-stream errors.
   Solution of the Euler equations at $t=4.0$  in a fully periodic domain, initialized with free-stream $(\rho,v_1,v_2,v_3,p)\!=\!(0.7,0.2,0.3,-0.4,1.0)$. }
   \label{fig:NonConformingPressFluc}
\end{figure}

Curved elements with non-conforming faces, however, can introduce geometrically induced errors, including free-stream preservation errors, that are not seen with straight-sided elements. Such geometrically induced errors do not have to be small. Fig. \ref{fig:NonConformingPressFluc} shows a simple mesh on which the pressure errors are 0.2\% using a fourth order interpolation polynomial to approximate both the solution and the geometry.

In this paper we show that free-stream preservation errors are generally introduced in the DGSEM approximation on non-conforming curved elements, as depicted in Fig.~\ref{fig:NonConformingPressFluc}. We then derive sufficient conditions to ensure free-stream preservation on non-conforming meshes where the elements are non-conforming due to subdivision of a parent element, e.g., Fig. \ref{fig:ConformingDiagrams}.  Since the conditions are purely geometrical, it is likely that other high order methods share the same issues, \modColor{e.g. high order flux reconstruction (FR) methods are similar to DG \cite{Allaneau:2011pz},\cite{Mengaldo2016}.}

\section{The Spectral Element Mesh}
 
\begin{figure}[tbp] 
   {\includegraphics[width=0.49\textwidth]{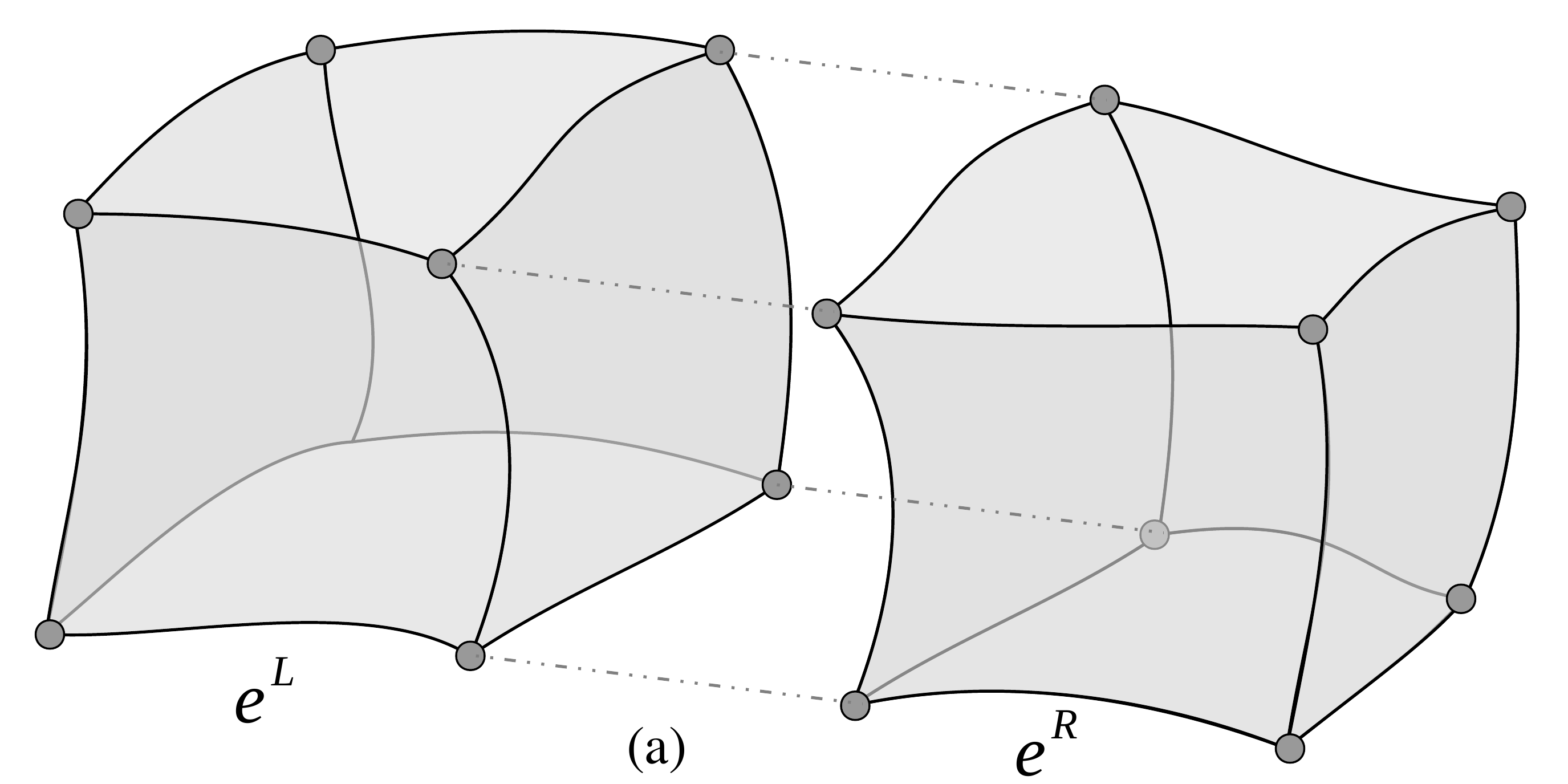} }
   {\includegraphics[width=0.49\textwidth]{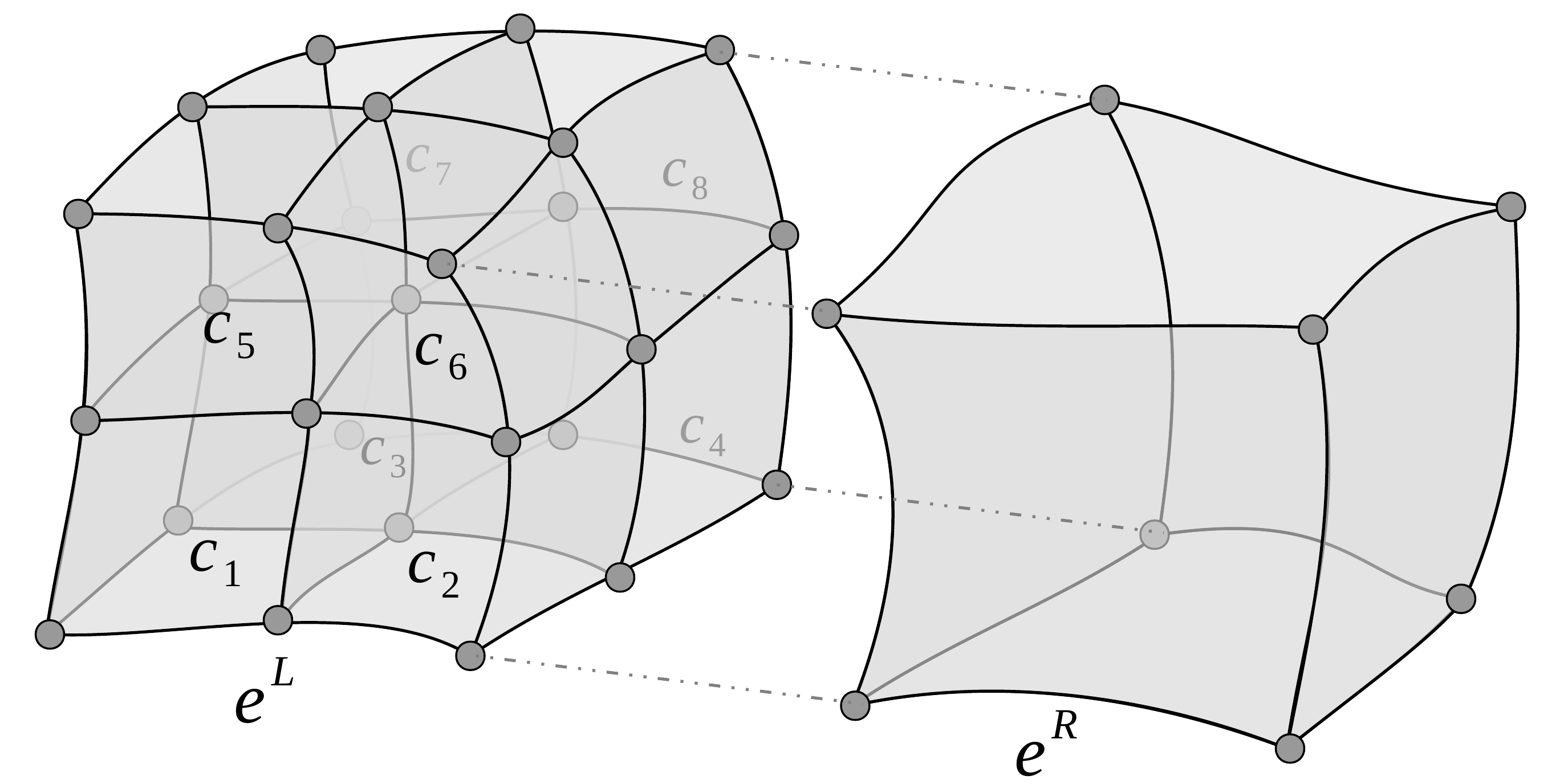} }
   \caption{Schematics of conforming elements (a) and geometrically non-conforming elements (b) where the left element $e^{L}$ is subdivided into eight child elements $c_{k}$. In general, the element faces will be curved. }
   \label{fig:ConformingDiagrams}
\end{figure}

The mesh is an unstructured subdivision of a domain $\Omega$ into non-overlapping hexahedral elements. 
 If all elements in a mesh share an entire edge, face, or only a corner, then the mesh is geometrically conforming (Fig. \ref{fig:ConformingDiagrams}a). Unstructured hexahedral element meshes are difficult to construct, but have desirable properties \cite{Kovalev:2005qv}. However, local refinement in a hexahedral mesh is most easily accomplished by introducing non-conforming elements (Fig. \ref{fig:ConformingDiagrams}b), where two elements might share only part of a face or edge. For our purposes here, it is sufficient to consider when a neighboring element $e^{L}$ is subdivided into eight child elements, $c_{i},\; i=1,2,\ldots,8$.

\begin{figure}[htbp] 
   \centering
   \includegraphics[trim=0 40 0 35,clip,width=3in]{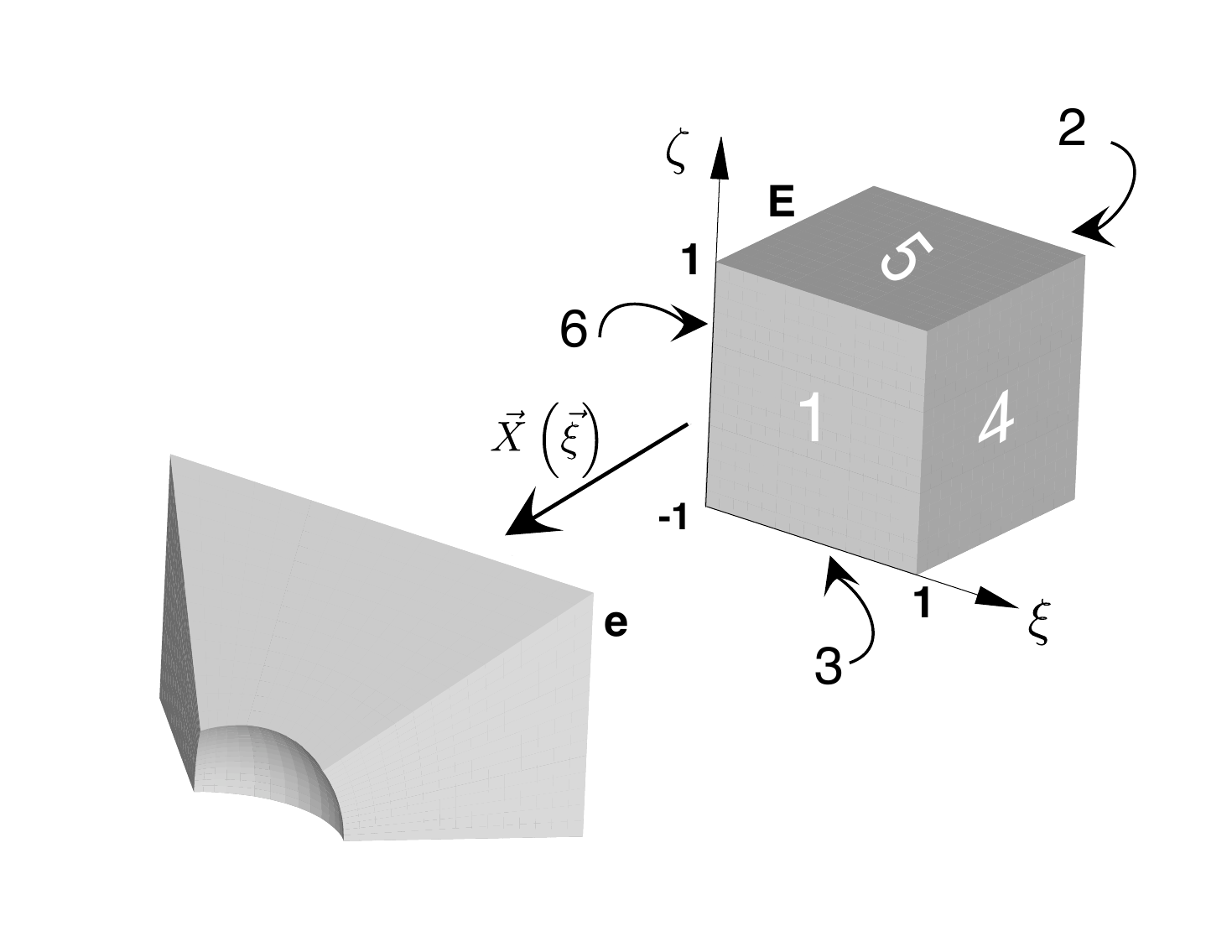} 
   \caption{The reference element with its face orientation and its mapping to an element in the mesh.}
   \label{fig:ElementMapping}
\end{figure}

To be used in the approximation, each element is mapped to a reference element $E=[-1,1]^{3}$. Precisely, one maps
 from the reference coordinates
 \begin{equation}
 \spacevec \xi = \sum\limits_{i = 1}^3 {{\xi^{i}}{{\hat \xi}^{i}}}  = \left( {\xi^{1}},{\xi^{2}},{\xi^{3}} \right) = (\xi,\eta,\zeta) =\xi\hat \xi + \eta\hat \eta + \zeta\hat \zeta\in [-1,1]^{3}
 \end{equation} 
to the element by a transformation  $\spacevec x = \spacevec X( \spacevec \xi\, )$, where $\spacevec x = (x,y,z) = \left(x_{1},x_{2},x_{3}\right) = \sum\limits_{n = 1}^3 {{x_n}{{\hat x}_n}} = x\hat x + y\hat y + z\hat z$. Fig. \ref{fig:ElementMapping} shows the geometry of the reference element and the mapping to an element.

For accurate approximation of physical boundaries, the mappings will generally not be tri-linear. Instead, the faces of $e$ along physical boundaries will be curved surfaces. Furthermore, to handle highly anisotropic meshes like those in boundary layer computations, curved elements will extend into the domain \cite{Persson:2009mz},\cite{Hindenlang:2014gl}. Thus, $\spacevec X$ will, in general, not be a linear function in each of its arguments. Typically, the boundaries are approximated by high order polynomial interpolants \cite{Patera:1984:JCP}. From those boundary interpolants, the mapping $\spacevec X$ is often a linear transfinite map \cite{Gordon&Hall1973a} that matches the element faces, although other types of maps might be useful \cite{ISI:000306588600006},\cite{Lovgren:2009zr}. 
 
To describe curved elements bounded by faces that are defined by polynomial interpolants, we introduce the polynomial interpolation operator. Let $\mathbb{P}^{N}=\mathbb{P}_{\xi}^{N}\times\mathbb{P}_{\eta}^{N}\times\mathbb{P}_{\zeta}^{N}$ be the tensor product space of polynomials of degree less than or equal to $N$ in each space direction. For some function $\statevec{v}\left(\spacevec\xi\right)$ defined on the reference element, the interpolant of $\statevec{v}$ through the tensor product of a one-dimensional distribution of nodes is
\begin{equation}
{\mathbb{I}^N}\left(\statevec{v}\left( \spacevec \xi \,\right)\right) = \sum\limits_{j,k,l = 0}^N {{\statevec{v}_{jkl}}{\ell _j}\left( \xi  \right){\ell _k}\left( \eta  \right){\ell _l}\left( \zeta  \right)} \in \mathbb{P}^{N},\end{equation}
where the $\ell_{j}$ is the Lagrange interpolating polynomial and ${\statevec{v}_{jkl}}$ is the value of $\statevec{v}$ at those nodes. The nodes will generally be the same as those used in the approximation of the solution described below in Sec. \ref{sec:DGSEM}, namely the Legendre Gauss or Gauss-Lobatto points.

 We then assume that the transformation for element $e^{m}$ is $\spacevec X^{m} \in\mathbb{P}^{N_{g}}\subset\mathbb{P}^{N}$, defined so that it matches the six faces , $\spacevec \Gamma_{s}\in\mathbb{P}^{N_{g}},\; s = 1, 2, \ldots,6$, i.e. $\spacevec X\left( -1, \eta, \zeta \right) = \spacevec \Gamma_{6}\left(\eta,\zeta\right)$, etc. It is important to remember that since $\mathbb{P}^{N_{g}}\subset\mathbb{P}^{N}$ for $N_{g}\le N$, the $\spacevec \Gamma_{s}$ can be of different polynomial orders $\left(N_{g}\right)_{s}\le N$ and still $\spacevec X\in\mathbb{P}^{N}$.

The transformations from the reference element to physical space must create a mesh that is \emph{watertight}, meaning that there are no gaps. The watertight condition constrains the mapping functions at faces of adjacent elements. For geometrically conforming meshes, the condition simply requires that the transformation along a face from the left and right elements be the same. 

Enforcing the watertight condition for non-conforming elements requires care. If a single analytic representation of the element faces is used, and if the transformation used in the interior matches that representation then the faces will match. However, since element faces will more likely be represented by high order polynomials, it is necessary to ensure that the faces represent the same polynomial. 

\begin{figure}[tbp] 
   \centering
   \includegraphics[width=3in]{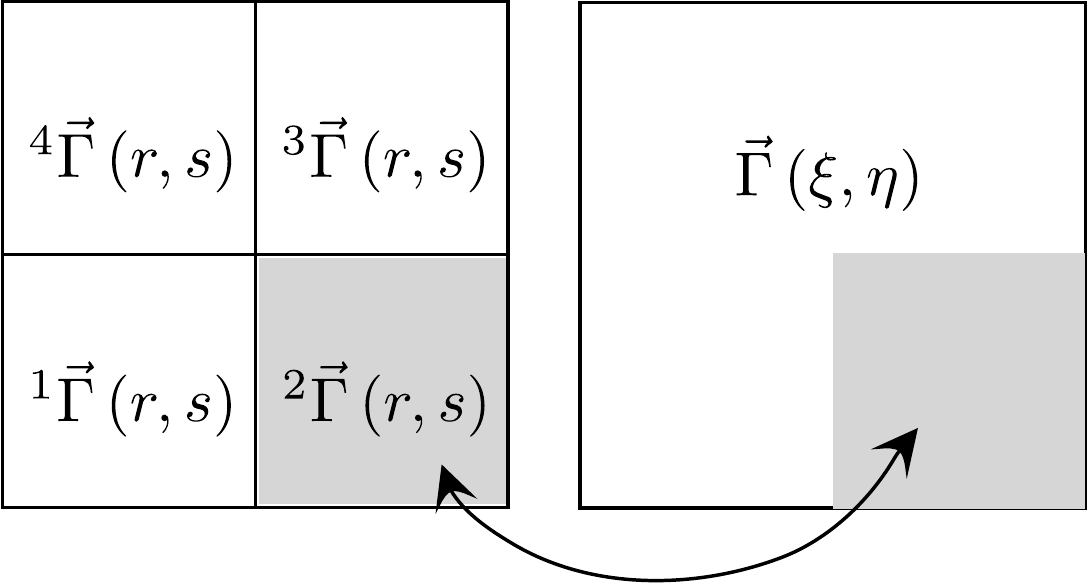} 
   \caption{Schematic representation of the coincident faces of the element $e^{L}$ and $e^{R}$ as shown in Fig. \ref{fig:ConformingDiagrams}b. We call the face on the right the \emph{parent face} and the faces on the left the \emph{child faces}.}
   \label{fig:SubdivisionMappings.pdf}
\end{figure}

As an example of how the face mappings must be defined, we consider the faces of a four-to-one subdivision as represented in Fig. \ref{fig:SubdivisionMappings.pdf}. The square on the right represents one of the six faces of an element. The left represents a matching face from a neighbor element that has been subdivided into eight child elements. For the mesh to be watertight, each child element face must match its portion on the right face. 

The face on the right of Fig. \ref{fig:SubdivisionMappings.pdf} is described by a single polynomial space, $\spacevec\Gamma\left(\xi,\eta\right)\in\mathbb{P}^{N_{g}}$. The four patches on the left describe a piecewise polynomial space. The piecewise space is the bigger space since one can represent the polynomial $\spacevec\Gamma\left(\xi,\eta\right)$ exactly by the four polynomials on the left, but not vice-versa. The piecewise polynomial space includes patches that have slope discontinuities at the intersections, which the global polynomial, $\spacevec\Gamma\left(\xi,\eta\right)$, cannot have. 

Since the full polynomial space is a subspace of the piecewise space, the four face patches on the left of Fig. \ref{fig:SubdivisionMappings.pdf} must be computed from the 
face on the right to ensure that the mesh is watertight:
\begin{equation}
{}^2\spacevec \Gamma \left( {r,s} \right) = \sum\limits_{i,j = 0}^{N_{g}} {{}^2{{\spacevec \Gamma }_{ij}}{\ell _i}\left( r \right){\ell _j}\left( s \right)} ,
\end{equation}
where
\begin{equation}
{}^2{\spacevec \Gamma _{ij}} = \spacevec \Gamma \left( {\frac{{{r_i} + 1}}{2},\frac{{{s_j} - 1}}{2}} \right),
\end{equation}
so that 
\begin{equation}{}^2\spacevec \Gamma \left( {r,s} \right) = \spacevec \Gamma \left( {\frac{{r + 1}}{2},\frac{{s - 1}}{2}} \right),
\label{eq:AffineMap}
\end{equation}
etc. for ${}^k\spacevec \Gamma,\;k=1,3,4$.

It is important to note that to represent the full space on the piecewise space exactly as a polynomial of degree $N_{g}$, the transformations of ${}^k\spacevec \Gamma, \;k=1,2,3,4$ onto their portions of $\spacevec \Gamma$, e.g. as seen in \ref{eq:AffineMap}, must be \emph{affine} to ensure that the polynomial order does not increase. Therefore, we will assume in what follows that, $\partial r/\partial \xi = \text{const}$ and $\partial s/\partial \eta = \text{const}$. From a practical point of view, this is the simplest way to subdivide elements anyway. (The constant does not have to be one half or the same for each coordinate direction.) 

For a watertight mesh: 
\begin{itemize}
\item Subdivide an element into children so that the transformation of each child to the parent is affine.
\item The child faces shared with their neighbor face must represent the same polynomial as on the neighbor face.
\item A three-dimensional mesh must also take the edge connectivity into account, so that edges on a finer level inherit the mapping of the lowest level on that edge.
\end{itemize}

\section{Coordinate Transformations, Metric Identities and Constant State Preservation}


 We next review the geometric conditions that are satisfied by the transformations and how they are related to free-stream preservation. From the transformation we can define two coordinate bases. The covariant basis
 \begin{equation}
\spacevec a_i  = \frac{{\partial \spacevec X}}
{{\partial \xi ^i }} = X_{\xi^{i}}\hat x + Y_{\xi^{i}}\hat y+ Z_{\xi^{i}}\hat z \quad i = 1,2,3
\end{equation}
 points along coordinate lines and can be derived directly from the transformation. The contravariant basis is normal to the coordinate lines,
\begin{equation}
{{\spacevec a}^i}= \nabla {\xi ^i} = \sum\limits_{n = 1}^3 {a_n^i{{\hat x}_n}}   = \frac{1}{\mathcal{J}}{{\spacevec a}_j} \times {{\spacevec a}_k}\quad i = 1,2,3\quad (i,j,k)\;cyclic,
\label{eq:ContravariantBasis}
\end{equation}
 where $\mathcal{J}=\spacevec a_{1}\cdot\left(\spacevec a_{2} \times \spacevec a_{3}\right)$ is the Jacobian of the transformation.
 It was shown in \cite{kopriva2006} that
 \begin{equation}
 \begin{split}
 \mathcal{J}a_n^i &=  - {{\hat x}_i} \cdot {\nabla _\xi } \times \left( {{X_l}{\nabla _\xi }{X_m}} \right) \\&= +{{\hat x}_i} \cdot {\nabla _\xi } \times \left( {{X_m}{\nabla _\xi }{X_l}} \right) \\&= \frac{1}{2}\;{{\hat x}_i} \cdot {\nabla _\xi } \times \left[ {{X_m}{\nabla _\xi }{X_l} - {X_l}{\nabla _\xi }{X_m}} \right]\quad i = 1,2,3;\;n = 1,2,3;\;(n,m,l)\;cyclic.
 \end{split}
 \label{eq:CurlForms}
 \end{equation}
  Written out explicitly in the first curl form, the metric terms are
 \begin{equation}
\begin{array}{l}
 \mathcal{J}\spacevec a^1  = \left[ {\left( {Y_\eta  Z} \right)_\zeta   - \left( {Y_\zeta  Z} \right)_\eta  } \right]\hat x + \left[ {\left( {Z_\eta  X} \right)_\zeta   - \left( {Z_\zeta  X} \right)_\eta  } \right]\hat y + \left[ {\left( {X_\eta  Y} \right)_\zeta   - \left( {X_\zeta  Y} \right)_\eta  } \right]\hat z \\ \\
 \mathcal{J}\spacevec a^2  = \left[ {\left( {Y_\zeta  Z} \right)_\xi   - \left( {Y_\xi  Z} \right)_\zeta  } \right]\hat x + \left[ {\left( {Z_\zeta  X} \right)_\xi   - \left( {Z_\xi  X} \right)_\zeta  } \right]\hat y + \left[ {\left( {X_\zeta  Y} \right)_\xi   - \left( {X_\xi  Y} \right)_\zeta  } \right]\hat z \\ \\
 \mathcal{J}\spacevec a^3  = \left[ {\left( {Y_\xi  Z} \right)_\eta   - \left( {Y_\eta  Z} \right)_\xi  } \right]\hat x + \left[ {\left( {Z_\xi  X} \right)_\eta   - \left( {Z_\eta  X} \right)_\xi  } \right]\hat y + \left[ {\left( {X_\xi  Y} \right)_\eta   - \left( {X_\eta  Y} \right)_\xi  } \right]\hat z \\ 
 \end{array}.
 \label{eq:JAiWrittenOut}
\end{equation}
Therefore, it follows that
 \begin{equation}
 \begin{split}
\sum\limits_{i = 1}^3 {\frac{{\partial \left( {\mathcal{J}a_n^i } \right)}}{{\partial \xi ^i }}}  &=  - \sum\limits_{i = 1}^3 {\frac{{\partial \left( { \hat x_i  \cdot \nabla _\xi   \times \left( {X_l \nabla _\xi  X_m } \right)} \right)}}{{\partial \xi ^i }}}  \\&=  - \nabla _\xi   \cdot \left( {\nabla _\xi   \times \left( {X_l \nabla _\xi  X_m } \right)} \right) = 0,\quad n = 1,2,3;\;(n,m,l)\;cyclic,
\end{split}
\end{equation}
giving us the so-called metric identities in vector form
\begin{equation}
\nabla_{\xi}  \cdot \left( {\mathcal{J}{{\spacevec a}^i}} \right) = 0.
\end{equation}

\begin{rem}
We note four important observations about the contravariant basis vectors: 
\begin{itemize}
\item Eq. (\ref{eq:ContravariantBasis}) says that at the element faces $\mathcal{J}\spacevec a^{i}$ is in the face normal direction. 
\item Eq.  (\ref{eq:JAiWrittenOut}) makes it easy to see that the derivatives that define the normal directions are in the face tangential directions. 
\item Since the transformation $\spacevec X$ must match the element faces at the faces of the reference element, it follows that at the faces the metric terms are computed solely from the derivatives of the parent face functions, $\spacevec \Gamma_{s}$.
\item With the affine relationship between the local coordinate on the child faces, $\spacevec \xi$, and the local coordinate on the full neighbor face, $\spacevec r$, the $\mathcal{J}\spacevec a^{i}$ computed on the child faces and the full neighbor face differ by a constant factor, $\beta$, i.e.
\begin{equation}
{\left. {\mathcal{J}{{\spacevec a}^i}} \right|_\text{child face}} = \frac{{\partial \spacevec X}}{{\partial {r ^j}}} \times \frac{{\partial \spacevec X}}{{\partial {r ^k}}} = \left( {\frac{{\partial {\xi^j}}}{{\partial {r ^j}}}\frac{{\partial {\xi^k}}}{{\partial {r ^k}}}} \right)\frac{{\partial \spacevec X}}{{\partial {r^j}}} \times \frac{{\partial \spacevec X}}{{\partial {r^k}}} \equiv \beta \cdot {\left. {\frac{{\partial \spacevec X}}{{\partial {\xi^j}}} \times \frac{{\partial \spacevec X}}{{\partial {\xi^k}}}} \right|_\text{parent face}}.\end{equation}
\end{itemize}
\label{rem:MetricFacts}
\end{rem}

The conservation law remains a conservation law in the reference domain coordinates,
\begin{equation}
\frac{{\partial \left( {\mathcal{J}{\statevec{u}}} \right)}}{{\partial t}} + {\nabla _\xi } \cdot \bigcontravec{f} = 0,
\label{eq:tranformedconsLaw}
\end{equation}
where $\bigcontravec{f}$ is the contravariant flux with components
\begin{equation}
 {\contrastatevec{f}}^{i} = \mathcal{J}{{\spacevec a}^i} \cdot \bigstatevec{f}  = \sum\limits_{n = 1}^3 {\mathcal{J}a_n^i{{\statevec{f}}_n}\hat x_{n}}.\end{equation}
Then if the covariant flux of a constant is constant, i.e., $\bigstatevec{f}\left({\statevec{c}}\right) = \bigstatevec{C}$,
\begin{equation}{\nabla _\xi } \cdot \bigcontravec{f} = {\nabla _\xi } \cdot \bigcontravec{C} = \sum\limits_{i = 1}^3 {\frac{{\partial \left( {\mathcal{J}{{\spacevec a}^i} \cdot \bigstatevec{C}} \right)}}{{\partial {\xi ^i}}}}  = \sum\limits_{i = 1}^3 {\frac{{\partial \left( {\mathcal{J}{{\spacevec a}^i}} \right)}}{{\partial {\xi ^i}}}}  \cdot \bigstatevec{C} = 0.
\end{equation}
Similarly, assuming the divergence of the flux of the constant state is zero implies the metric identities. 
In other words, it is both necessary and sufficient that the metric identities be satisfied for the transformed equations to also be constant state preserving.

\section{Discontinuous Galerkin Spectral Element Approximations}\label{sec:DGSEM}
We now review the construction of discontinuous Galerkin spectral element methods (DGSEM) for conforming and non-conforming hexahedral elements. Details of the approximation can be found in \cite{Kopriva:2009nx} for the standard approximation. The non-conforming approximation is described in \cite{Koprivaetal2000}.

\subsection{The Spatial Approximation within an Element}\label{sec:SpatialApproximation}

The DGSEM approximates a weak form of the equations \eqref{eq:tranformedconsLaw} on each element. The conservation law is multiplied by a test function $\statevecg{\phi} \in \mathbb{L}^{2}(E)$, where 
\begin{equation}{\mathbb{L}^2}(E) = \left\{ {{\statevec{u}}:\iprod {{\statevec{u}},{\statevec{u}}} = \int_E {{{\statevec{u}}^T}{\statevec{u}}d\spacevec \xi }  = {{\left\| {\statevec{u}} \right\|}^2} < \infty } \right\},
\end{equation}
and integrated over the reference domain. In inner product notation, the weak form is
\begin{equation}
\iprod {\mathcal{J}{{\statevec{u}}_t},\statevecg \phi } + \iprod{{\nabla _\xi } \cdot \bigcontravec{f},\statevecg \phi }  = 0.
\end{equation}
We then integrate the second inner product by parts to separate the volume from surface contributions and get the weak form to be approximated
\begin{equation}\iprod {\mathcal{J}{{\statevec{u}}_t},\statevecg \phi }+ {\left. {{\bigcontravec{f}^T}\statevecg \phi } \right|_{\partial E}} - \iprod {\bigcontravec{f},{\nabla _\xi }\statevecg \phi }= 0.
\label{eq:WeakformToApproximate}
\end{equation}
Here we have introduced the shorthand notation
\begin{equation}
{\left. {{\bigcontravec{f}^T}\statevecg \phi } \right|_{\partial E}} \equiv \sum\limits_{s = 1}^6 {\left\{ {\int_{\text{face}^{s}} {{{\left( {\bigcontravec{f} \cdot {{\hat n}_{\xi}^s}} \right)}^T}\statevecg \phi \,\dS^s} } \right\}} ,
\label{eq:boundaryTermsToApproximate}
\end{equation}
where $\hat n_{\xi}^{s}$ is the outward normal in the reference space on face $s$. The normal component of the contravariant flux, $\bigcontravec{f}\cdot\hat n^{s}_{\xi}$, is proportional to the normal flux
\begin{equation}\bigcontravec{f} \cdot \hat n_\xi ^s = \sum\limits_{n = 1}^3 {Ja_n^s{{\statevec{f}}_n}}  = \left| {\nabla {\xi ^s}} \right|\bigstatevec{f}  \cdot {{\hat n}_x},\end{equation}
where $|\cdot |$ is the Euclidean magnitude of the vector.


The discontinuous Galerkin SEM approximates the solutions and fluxes in (\ref{eq:WeakformToApproximate}) by polynomials of degree $N$ written in Lagrange form with nodes at the quadrature points, takes $\statevecg \phi$ to be a polynomial of degree $N$, replaces integrals by Gauss or Gauss-Lobatto quadratures and replaces the boundary flux in (\ref{eq:boundaryTermsToApproximate}) by a numerical flux (Riemann solver) that couples neighboring elements. 
We therefore approximate 
\begin{equation}\begin{gathered}
 {\statevec{u}}  \approx {\statevec{U}} = \sum\limits_{j,k,l = 0}^N {{{\statevec{U}}_{jkl}}{\ell _j}\left( \xi  \right){\ell _k}\left( \eta  \right){\ell _l}\left( \zeta  \right)}  \in {\mathbb{P}^N} \hfill \\
   \mathcal{J} \approx J = \sum\limits_{j,k,l = 0}^N {{J_{jkl}}{\ell _j}\left( \xi  \right){\ell _k}\left( \eta  \right){\ell _l}\left( \zeta  \right)}  \in {\mathbb{P}^N} \hfill \\
{{\contrastatevec{f}}^i} \approx {{\contrastatevec{F}}^i} = \sum\limits_{j,k,l = 0}^N {{\contrastatevec{f}}_{jkl}^i{\ell _j}\left( \xi  \right){\ell _k}\left( \eta  \right){\ell _l}\left( \zeta  \right)}  \in {\mathbb{P}^N}, \hfill \\
\end{gathered}
\end{equation} 
where
\begin{equation}
{\contrastatevec{f}}_{jkl}^i = \sum\limits_{n = 0}^3 {{{\statevec{f}}_n}\left( {{{\statevec{U}}_{jkl}}} \right){{\left( {J\spacevec a_n^i} \right)}_{jkl}}}  = {\bigstatevec f _{jkl}} \cdot {\left( {J{{\spacevec a}^i}} \right)_{jkl}}.\end{equation}
Note that the interpolant of the contravariant flux is equivalently written as ${ {\contrastatevec{F}}^i} = {\mathbb{I}^N}\left( {{\mathbb{I}^N}\left( \bigstatevec{f} \right)\cdot{\mathbb{I}^N}\left( {J{{\spacevec a}^i}} \right)} \right)$.

We also define the discrete inner product of two polynomials, $\statevec{ V},\statevec{ W}\in\mathbb{P}^{N}$
as
\begin{equation}
{\iprod {\statevec{ V},\statevec{ W}} _N} = \sum\limits_{j,k,l = 0}^N {\statevec{ V}_{jkl}^T{{\statevec{ W}}_{jkl}}{\omega_j}{\omega_k}{\omega_l}}  \equiv \sum\limits_{j,k,l = 0}^N {\statevec{ V}_{jkl}^T{{\statevec{ W}}_{jkl}}{\omega_{jkl}}},
 \end{equation}
where the singly subscripted $\omega$'s are the one-dimensional Legendre-Gauss(-Lobatto) quadrature weights and the triply subscripted is the product of the three. We use a similar notation for integrals, where we add a subscript $N$ to denote quadrature, e.g.
\begin{equation}
\int_{E,N} {V \mathrm{d}\spacevec \xi }  \equiv \sum\limits_{j,k,l = 0}^N {{V_{jkl}}{\omega_{jkl}}} .
\end{equation}

With this notation, the formal statement of the DGSEM semi-discrete approximation on an element becomes
\begin{equation}
{\iprod{ {J{{\statevec{U}}_t},\statevecg \phi } }_N} + {\left. {{{ {\contrastatevec{F}}}}^{*,T}\statevecg \phi } \right|_{\partial E,N}} - {\iprod {\bigcontravec{F},{\nabla _\xi }\statevecg \phi } _{N}} = 0,
\label{eq:FirstWeakForm}
\end{equation}
with the contravariant numerical flux
\begin{equation}
{\contrastatevec{F}}^* = {{\statevec{F}}^*}\left( {{{\statevec{U}}^L},{{\statevec{U}}^R};J{{\spacevec a}^s}} \right)
\label{eq:FaceContravariantFlux}
\end{equation}
derived from the usual normal Riemann flux, $\statevec{F}^{*}\left( {{{\statevec{U}}^L},{{\statevec{U}}^R};{{\hat n}_x}} \right)$. \modColor{If one replaces the volume and surface flux quadratures with $M>N$, typically called ``overintegration'', one provides the chance to increase the precision of the quadrature to eliminate or reduce aliasing errors associated with the representation of the contravariant flux as a polynomial of degree $N$.} 

The discrete weak form (\ref{eq:FirstWeakForm}) can be written in an alternate penalty form. The use of Gauss-type quadratures leads to a summation by parts formula \cite{gassner2010}, precisely,
\begin{equation}{\iprod { {\bigcontravec{F}},{\nabla _\xi }\statevecg \phi } _{N}} = {\left. {{\left({ {\bigcontravec{F}}\cdot\hat n}\right)^T}\statevecg \phi } \right|_{\partial E,N}} - {\iprod {{\nabla _\xi } \cdot {\bigcontravec{F}},\statevecg \phi } _{N}}.
\end{equation}
Therefore, we can also write the algebraically equivalent form
\begin{equation}{\iprod {J{{\statevec{U}}_t},\statevecg \phi } _N} + {\left. {{{\left\{ {{{ {\contrastatevec{F}}}^*} -  {\bigcontravec{F}}\cdot \hat n} \right\}}^T}\statevecg \phi } \right|_{\partial E,N}} + {\iprod {{\nabla _\xi } \cdot  {\bigcontravec{F}},\statevecg \phi } _{N}} = 0.
\label{eq:formII}
\end{equation}

The equivalent forms \eqref{eq:FirstWeakForm} and \eqref{eq:formII} are not necessarily stable \cite{Kopriva:2017yg}, even for conforming approximations. For a linear system of equations where the flux components are of the form $\statevec F_{i} = \mmatrix A_{i}\statevec U$, $\mmatrix A_{i}= \text{const}$, a strong form skew-symmetric approximation that is stable for conforming approximations is \cite{Kopriva2016274}
 \begin{equation}
 \begin{split}
 \quad {\iprod {J\statevec{U}_{t},\statevecg \phi }}_N   
 &+ \frac{1}{2}\left\{ 
   {{\iprod {\nabla  \cdot \bigcontravec F \left( \statevec U \right), \statevecg \phi }}_N} 
 + {{\iprod {{\mathbb{I}^N}\left( {\tilde {\mmatrix A}} \right) \cdot \nabla\statevec{U}, \statevecg \phi  }}_N} 
 + {{\iprod {\nabla \cdot \left({\mathbb{I}^N}\left( {\tilde {\mmatrix A}} \right)\right){\statevec{U}},\statevecg \phi }}_N} \right\} \\
 &+ \int_{\partial E,N} {\left\{ \contrastatevec{F} ^* - \bigcontravec{F} \cdot \hat n \right\}^T}\statevecg \phi \dS = 0.
 \end{split}
 \label{eq:linearSplitDGSEM}
 \end{equation}
\modColor{ Alternatively, overintegration of the linear problem leads to a stable approximation \cite{Kopriva2018},\cite{Mengaldo201556} on conforming meshes.}
 
For the nonlinear Euler equations, discontinuous Galerkin approximation that is entropy stable for conforming approximations is \cite{Gassner2018}
\begin{equation}
\label{eq:disc_nse1}
\iprod{J\,\statevec{U}_t,\statevecg{\phi}}_N + \iprod{\spacevec{\mathbb{D}} (\bigcontravec{F} )^{\ec},\statevecg{\phi}}_N +\int\limits_{\partial E,N} \statevecg{\phi}^T\left(\contrastatevec{F}^{\ec,*}-\bigcontravec{F}\cdot\hat{n}\right)\,\dS =0. 
\end{equation}
The $\contrastatevec{F}^{\ec,*}$ is an entropy conservative numerical flux. Using the definition 
\begin{equation}
\average{J\spacevec{a}^{\,1}}_{(i,m)jk}\equiv \frac{1}{2}\left\{J\spacevec{a}^{\,1}_{ijk}+J\spacevec{a}^{\,1}_{mjk}\right\},
\end{equation}
etc.,
\begin{equation}
\label{eq:entropy-cons_volint}
\begin{split}
\spacevec{\mathbb{D}} (\bigcontravec{F} )^{\ec}_{ijk}\equiv
&\quad 2\sum_{m=0}^N D_{im}\,\left(\bigstatevec{F}^{\ec}(\statevec U_{ijk}, \statevec U_{mjk})\cdot\average{J\spacevec{a}^{\,1}}_{(i,m)jk}\right)\\[0.1cm]
&+     2\sum_{m=0}^N D_{jm}\,\left(\bigstatevec{F}^{\ec}(\statevec U_{ijk}, \statevec U_{imk})\cdot\average{J\spacevec{a}^{\,2}}_{i(j,m)k}\right)\\[0.1cm]
&+     2\sum_{m=0}^N D_{km}\,\left(\bigstatevec{F}^{\ec}(\statevec U_{ijk}, \statevec U_{ijm})\cdot\average{J\spacevec{a}^{\,3}}_{ij(k,m)}\right),
\end{split}
\end{equation}
where $\bigstatevec{F}^{\ec}(\statevec U_{ijk},\statevec U_{mjk})$ is a two-point entropy conserving flux. \modColor{For the purposes of this paper, the entropy conservative numerical flux and the two-point numerical flux are \emph{consistent}. They are defined so that $\contrastatevec{F}^{\ec,*}\left(\statevec U, \statevec U\right) = \bigstatevec f(\statevec U)\cdot\hat n$ and $ \bigstatevec{F}^{\ec}(\statevec U,\statevec U) = \bigstatevec{f}(\statevec U)$. For further properties, see \cite{Gassner2018}, \cite{Gassner:2016ye}}

\section{The Discrete Metric Identities and Free-Stream Preservation}
We are concerned with what happens when the solution is constant, i.e., $\statevec{U}=\statevec{c}$ so that the flux $\bigstatevec{f}=\bigstatevec{C} $ is constant. The approximation should imply that $\statevec{U}_{t}=0$. 

To derive the conditions needed for the approximate solution to also remain constant, we first expand the divergence in the volume term of the second form of the approximation, (\ref{eq:formII})
\begin{equation}{\nabla _\xi } \cdot \bigcontravec{F} = \sum\limits_{i = 1}^3 {\frac{{\partial {{ {\contrastatevec{F}}}^i}}}{{\partial {\xi ^i}}}}  = \sum\limits_{n = 1}^3 {\sum\limits_{i = 1}^3 {\frac{\partial }{{\partial {\xi ^i}}}} \left\{ {{\mathbb{I}^N}\left( {{\mathbb{I}^N}\left( {{{\statevec{f}}_n}} \right){\mathbb{I}^N}\left( {Ja_n^i} \right)} \right)} \right\}}. \end{equation}
With $\statevec{f}_{n}=\statevec{C}_{n} = const$
\begin{equation}{\nabla _\xi } \cdot \bigcontravec{F}\left( {\statevec{c}} \right) = \sum\limits_{n = 1}^3 {{{\statevec{C}}_n}\sum\limits_{i = 1}^3 {\frac{\partial }{{\partial {\xi ^i}}}\left\{ {{\mathbb{I}^N}\left( {Ja_n^i} \right)} \right\}} }. \end{equation}
In the linear problem, \eqref{eq:linearSplitDGSEM}, 
\[
{\nabla _\xi } \cdot {\mathbb{I}^N}\left( {\tilde {\mmatrix A}} \right)\statevec c = \sum\limits_{n = 1}^3 {\sum\limits_{i = 1}^3 {\frac{\partial }{{\partial {\xi ^i}}}\left\{ {{\mathbb{I}^N}\left( {{\mmatrix A_n}Ja_n^i} \right)} \right\}} }\statevec c  = \sum\limits_{n = 1}^3 {{{\statevec{C}}_n}\sum\limits_{i = 1}^3 {\frac{\partial }{{\partial {\xi ^i}}}\left\{ {{\mathbb{I}^N}\left( {Ja_n^i} \right)} \right\}} },
\]
so
\[
\frac{1}{2}\left\{ {{{\left( {{\nabla  } \cdot \bigcontravec{F}\left( {\statevec{c}} \right),\statevecg \phi } \right)}_N} + {{\left( {{\mathbb{I}^N}\left( {\tilde {\mmatrix A}} \right) \cdot {\nabla  }{\statevec{c}},\statevecg \phi } \right)}_N} + {{\left( {{\nabla  } \cdot \left({\mathbb{I}^N}\left( {\tilde {\mmatrix A}} \right)\right){\statevec{c}},\statevecg \phi } \right)}_N}} \right\} = \sum\limits_{n = 1}^3 {{{\statevec{C}}_n}\sum\limits_{i = 1}^3 {\frac{\partial }{{\partial {\xi ^i}}}\left\{ {{\mathbb{I}^N}\left( {Ja_n^i} \right)} \right\}} }.
\]
Finally, since the discrete derivative of a constant is zero, and the two point flux is consistent with the original Euler flux,
\begin{equation}
\label{eq:entropy-cons_volintconst}
\begin{split}
\spacevec{\mathbb{D}} (\bigcontravec{F} )^{\ec}_{ijk}\equiv
&\quad 2\sum_{m=0}^N D_{im}\,\left(\bigstatevec{F}^{\ec}(\statevec c, \statevec c)\cdot\average{J\spacevec{a}^{\,1}}_{(i,m)jk}\right)\\[0.1cm]
&+     2\sum_{m=0}^N D_{jm}\,\left(\bigstatevec{F}^{\ec}(\statevec c, \statevec c)\cdot\average{J\spacevec{a}^{\,2}}_{i(j,m)k}\right)\\[0.1cm]
&+     2\sum_{m=0}^N D_{km}\,\left(\bigstatevec{F}^{\ec}(\statevec c, \statevec c)\cdot\average{J\spacevec{a}^{\,3}}_{ij(k,m)}\right)\\&
= \sum\limits_{n = 1}^3 {{{\statevec{C}}_n}\sum\limits_{l = 1}^3 {\frac{\partial }{{\partial {\xi ^l}}}\left\{ {{\mathbb{I}^N}\left( {Ja_n^l} \right)} \right\}} }_{ijk}.
\end{split}
\end{equation}

The surface terms for all three approximations become
\begin{equation}{\left. {{{\left\{ {{{\contrastatevec{F}}^*}\left( {\statevec{c}} \right) - \bigcontravec{F}\left( {\statevec{c}} \right)}\cdot\hat n \right\}}^T}\statevecg \phi } \right|_{\partial E,N}} = \sum\limits_{s = 1}^6 {\left\{ {\int_{\text{face}^{s},N} {\sum\limits_{n = 1}^3 {{\statevec{C}}_n^T\statevecg \phi \sum\limits_{i = 1}^3 {{{\hat \xi }^i}\left[ {{{\left( {{\mathbb{I}^N}\left( {Ja_n^i} \right)} \right)}^*} - \left( {{\mathbb{I}^N}\left( {Ja_n^i} \right)} \right)} \right]} }  \cdot {{\hat n}^s}\dS^s} } \right\}} \end{equation}where the $^*$ indicates the metric terms as represented on the face.

Since the vector $\bigstatevec{C}$ is arbitrary, it follows that for any of the approximations \eqref{eq:formII}, \eqref{eq:linearSplitDGSEM}, or \eqref{eq:disc_nse1}, a constant state is preserved if and only if the metric terms satisfy the \emph{Discrete Metric Identities}
\begin{equation}
\sum\limits_{s = 1}^6 {\left\{ {\int_{\text{face}^s,N} {\left( {\sum\limits_{i = 1}^3 {{{\hat \xi }^i}\left[ {{\mathbb{I}^N}\left( {{{\left( {J{{\spacevec a}^i}} \right)}^*}} \right) - \left( {{\mathbb{I}^N}\left( {J{{\spacevec a}^i}} \right)} \right)} \right]} } \right) \cdot {{\hat n}^s}\phi \dS^s} } \right\}}  - {\left( {\sum\limits_{i = 1}^3 {\frac{{\partial {\mathbb{I}^N}\left( {J{{\spacevec a}^i}} \right)}}{{\partial {\xi ^i}}},\phi } } \right)_N} = 0.
\label{eq:DMI}
\end{equation}

\begin{rem}
\modColor{To get the metric identities for an overintegrated approximation, the quadrature $N$ is replaced by $M>N$.}
\end{rem}
\begin{rem}
Note that \eqref{eq:DMI} is a purely geometric condition, independent of the original PDE being solved and the form \eqref{eq:formII}, \eqref{eq:linearSplitDGSEM}, or \eqref{eq:disc_nse1} of the approximation. The numerical flux must be computed in such a way that 
\eqref{eq:DMI} is satisfied for a free stream to be preserved. It is a consistency condition. Coupling the numerical fluxes in a stable manner is still an open problem for curved interfaces, especially for nonlinear problems, Cf. \cite{L.-Friedrichs:2017hm,Lucas-Friedrich:2018ng}. Nevertheless, a stable numerical flux procedure will not be free-stream preserving unless \eqref{eq:DMI} is satisfied. Conversely, satisfying \eqref{eq:DMI} does not guarantee stability.
\end{rem}

From a practical point of view, a constant state is preserved if the volume and surface terms vanish independently, that is, if
\begin{equation}
{\bf Condition \;(V):}  \quad{\left( {\sum\limits_{i = 1}^3 {\frac{{\partial {\mathbb{I}^N}\left( {J{{\spacevec a}^i}} \right)}}{{\partial {\xi ^i}}},\phi } } \right)_N} = 0
  \label{eq:VolumeMetricIdentity}
 \end{equation}
 and
 \begin{equation}
\int_{\text{face}^s,N} {\left( {\sum\limits_{i = 1}^3 {{{\hat \xi }^i}\left[ {{\mathbb{I}^N}\left( {{{\left( {J{{\spacevec a}^i}} \right)}^*}} \right) - \left( {{\mathbb{I}^N}\left( {J{{\spacevec a}^i}} \right)} \right)} \right]} } \right) \cdot {{\hat n}^s}\phi \dS^s}  = 0\quad s = 1,2, \ldots ,6
  \label{eq:SurfacesMetricIdentityI}
\end{equation}
hold.

Condition (V) determines how the metric terms $J\spacevec a^{i}$ must be computed within the volume of an element \cite{kopriva2006}. It states that the divergence of the interpolant of the metric terms must vanish. The condition can be enforced by using any one of the curl forms in (\ref{eq:CurlForms}). For instance, the metric terms can be approximated as
\begin{equation}
Ja_n^i =  - {\hat x_i} \cdot {\nabla _\xi } \times \left( {{\mathbb{I}^N}\left( {{X_l}{\nabla _\xi }{X_m}} \right)} \right) \in \mathbb{P}^{N},
\label{eq:PolynomialMetrics}
\end{equation}
e.g.
\begin{equation}
\begin{split}
J{{\spacevec a}^3} &= \left[ {{{\left( {{\mathbb{I}^N}\left( {{Y_\xi }Z} \right)} \right)}_\eta } - {{\left( {{\mathbb{I}^N}\left( {{Y_\eta }Z} \right)} \right)}_\xi }} \right]\hat x \\&+ \left[ {{{\left( {{\mathbb{I}^N}\left( {{Z_\xi }X} \right)} \right)}_\eta } - {{\left( {{\mathbb{I}^N}\left( {{Z_\eta }X} \right)} \right)}_\xi }} \right]\hat y \\&+ \left[ {{{\left( {{\mathbb{I}^N}\left( {{X_\xi }Y} \right)} \right)}_\eta } - {{\left( {{\mathbb{I}^N}\left( {{X_\eta }Y} \right)} \right)}_\xi }} \right]\hat z.
\end{split}
\end{equation}
Condition (V) is then satisfied because the $J\spacevec a^{i}$ are already polynomials of degree $N$ when computed this way, the interpolation in (\ref{eq:VolumeMetricIdentity}) has no effect, and the divergence is explicitly zero.  

Using one of the curl forms (\ref{eq:CurlForms}) is the most general way to enforce the Condition (V), (\ref{eq:VolumeMetricIdentity}).
Special cases for other approximate forms are discussed in \cite{kopriva2006}. Note that since $\spacevec X\in\mathbb{P}^{N_{g}}$, there is an aliasing error associated with the interpolation of the products in (\ref{eq:PolynomialMetrics}) if $N_{g}=N$.
Also note that the choice of the approximation in the volume fixes the surface values of ${{\mathbb{I}^N}\left( {J{{\spacevec a}^i}} \right)}$ in (\ref{eq:SurfacesMetricIdentityI}).

The surface integrals along a face $s$, whose normal is in the $\pm \hat\xi^{i}$ local coordinate direction, reduce to
\begin{equation}
{\bf Condition\; (F):}\quad\int_{\text{face},N} \left[ {{{\left( {{\mathbb{I}^N}\left( {J{{\spacevec a}^i}} \right)} \right)}^*} - \left( {{\mathbb{I}^N}\left( {J{{\spacevec a}^i}} \right)} \right)} \right]\phi \dS  =0.
\label{eq:SurfacemetricIdentityII}
\end{equation}
Therefore, the surface term vanishes if either the interpolant of the metric terms evaluated at the boundary matches those used to compute the numerical flux, or the difference is always orthogonal to the polynomial space on the element face. 

We can then state: \emph{The approximation is free-stream preserving if Condition (V) and Condition (F) hold}.

\subsection{Element Coupling:  Enforcing {\bf Condition (F)}}

The elements are coupled through the numerical flux, ${\contrastatevec{F}}^{*}$, created from the states on either side of the face. For conforming faces the nodes at which the solutions are represented are co-located so that the two states and face normal are unique, allowing one to compute a unique numerical flux at each point on a face. Recalling Rem. \ref{rem:MetricFacts}, the fact that conforming faces match means that the $\mathbb{I}^{N}\left(J\spacevec a^{i}\right)$ match from either side. It then makes sense to choose ${{\mathbb{I}^N}\left( {{{\left( {J{{\spacevec a}^i}} \right)}^*}} \right)} = \mathbb{I}^{N}\left(J\spacevec a^{i}\right)$ and Condition (F) is automatically satisfied.

When the faces are non-conforming as in Fig. \ref{fig:ConformingDiagrams}b then the nodes are not co-located, the solutions lie in different spaces, and a choice for ${{\mathbb{I}^N}\left( {{{\left( {J{{\spacevec a}^i}} \right)}^*}} \right)} $ must be made to ensure that Condition (F) is satisfied simultaneously on both elements that share the face.

A procedure is therefore needed to define the metric terms on the faces. One can either use the metric terms computed from $\spacevec \Gamma$ or from the subdivided child faces, ${}^{s}\spacevec \Gamma$. Again, the faces can represent a polynomial approximation from $\spacevec \Gamma$ but not vice-versa. 
For that reason, we choose $\left(J\spacevec a^{i}\right)^{*}$, like the face polynomial itself, to be the value computed from $\Gamma$, that is, we use the interpolation in the large element $e^{R}$
\begin{equation}
{\left( {Ja_n^i} \right)^*} = {\left. { - {{\hat x}_i} \cdot {\nabla _\xi } \times \left( {\mathbb{I}_{{e^R}}^N\left( {{\Gamma _l}{\nabla _\xi }{\Gamma _m}} \right)} \right)} \right|_\text{face}} \in {\mathbb{P}^N}
\label{eq:MortarmetrixApprox}
\end{equation}
to minimize the number of errors being introduced. The other curl forms of the metric terms can also be used, with appropriate modifications in what follows.

In general, Condition (F) will not hold for geometrically non-conforming elements. Using (\ref{eq:MortarmetrixApprox}) from which to compute the face value, the quadrature vanishes for the face $\Gamma$ because the projection from the faces preserves the polynomials. It does not vanish, for the faces ${}^{i}\Gamma$, however, unless 
\begin{equation}
{\left( {Ja_n^i} \right)^*} = {\left. { - {{\hat x}_i} \cdot {\nabla _\xi } \times \left( {\mathbb{I}_{{e^R}}^N\left( {{\Gamma _l}{\nabla _\xi }{\Gamma _m}} \right)} \right)} \right|_\text{face}} = {\left. { - \frac{1}{\beta}{{\hat x}_i} \cdot {\nabla _\xi } \times \left( {\mathbb{I}_{{c_k}}^N\left( {^{k}{\Gamma _l}{\nabla _\xi }{^{k}\Gamma _m}} \right)} \right)} \right|_\text{face}}
\end{equation} 
for the appropriate face ${}^{i}\Gamma$ associated with child element $c_{k}$ that are neighbors to $e^{R}$. Therefore, we have

\emph{Condition (F) holds if the tangential components satisfy}
\begin{equation}
\mathbb{I}_{{c_k}}^N\left( {^{k}{\Gamma_l}{\nabla _\xi }^{k}{\Gamma_m}} \right) =\beta\cdot\mathbb{I}_{{e^R}}^N\left( {{\Gamma_l}{\nabla _\xi }{\Gamma_m}} \right) 
\label{eq:GeneralCondition}
\end{equation}
\emph{at all faces shared by two elements.}

In general (\ref{eq:GeneralCondition}) will not hold if the metric terms are computed using the local interpolation $\mathbb{I}_{{c_k}}^N$ even if the faces match because the interpolation of the product of two polynomials onto elements of different size will differ (see the Appendix) unless, as noted above, the product ${{\Gamma_l}{\nabla _\xi }{\Gamma_m}}$ is already a polynomial of degree $N$. 

\begin{rem}
\modColor{As noted in Sec. \ref{sec:SpatialApproximation}, overintegration replaces the volume and surface integrals for the spatial terms with quadratures with $M>N$. With overintegration, Condition (V) still holds with $M\ge N$  for any interpolation of $X_{l}\nabla_{\xi}X_{m}$ in \eqref{eq:PolynomialMetrics}. Furthermore, one has the flexibility to choose a sufficient number of interpolation points so that with $N_{g}\le N$, 
\begin{equation}
\mathbb{I}_{{c_k}}^M\left( {^{k}{\Gamma_l}{\nabla _\xi }^{k}{\Gamma_m}} \right) =\beta\cdot\mathbb{I}_{{e^R}}^M\left( {{\Gamma_l}{\nabla _\xi }{\Gamma_m}} \right) 
=\beta\cdot\left( {{\Gamma_l}{\nabla _\xi }{\Gamma_m}} \right)
\label{eq:GeneralCondition2}
\end{equation}
The sufficient condition is therefore $M\ge 2N_{g}$.
}
\end{rem}

\section{Geometry Definitions for which a Free-Stream is Preserved}
\label{sec:catalog}
In this section we catalog element geometry definitions for which \eqref{eq:GeneralCondition} and Condition (F) are satisfied. We assume that the metric terms are computed so that Condition (V) is satisfied, e.g. with \eqref{eq:PolynomialMetrics}. Then a constant state is preserved for:
\begin{enumerate}
\item Conforming Elements\\
If two neighboring elements are conforming, the elements are watertight and the same order polynomial is used in each, then $\beta = 1$, $\mathbb{I}^{N}_{c_{1}} \equiv \mathbb{I}^{N}_{e^{L}} = \mathbb{I}^{N}_{e^{R}}$. Therefore, the interpolation projections on the left and right elements are the same so (\ref{eq:GeneralCondition}) holds. This was the situation demonstrated in \cite{kopriva2006}.
\item Elements with Bilinear Faces\\
As a corollary to the previous condition, if the element faces have straight edges, a bilinear mapping, and the polynomial order of the approximation is $N\ge 2$. Then (\ref{eq:GeneralCondition}) holds and the approximation is free-stream preserving.
\item Two-Dimensional Problems\\
In two space dimensions, the calculation of the metric terms reduces to
\[\begin{gathered}
  J{{\spacevec a}^1} = {Y_\eta }\hat x - {X_\eta }\hat y\quad  \in {\mathbb{P}^N} \hfill \\
  J{{\spacevec a}^2} =  - {Y_\xi }\hat x + {X_\xi }\hat y\; \in {\mathbb{P}^N} \hfill \\ 
\end{gathered}. \]
When the element mappings are polynomials of degree $N$, the metric terms are, too, and can be represented exactly on all faces (edges in 2D).
\item Two-Dimensional Extrusions\\
For two dimensional geometries extruded in the $z$ direction, $Z = \zeta \in \mathbb{P}_{\zeta}^{1}$. The other two components are tensor products of two one dimensional polynomials with $X=X\left(\xi,\eta\right)\in\mathbb{P}_{\xi}^{N}\times\mathbb{P}_{\eta}^{N}$ and $Y = Y\left(\xi,\eta\right)\in\mathbb{P}_{\xi}^{N}\times\mathbb{P}_{\eta}^{N}$. Then ${{X_l}{\nabla _\xi }{X_m}}\in \mathbb{P}^{N}$ and (\ref{eq:GeneralCondition}) holds. Therefore, approximations like those seen in \cite{Zhang2015147} should be free-stream preserving.
\item Parent element metrics\\
For those meshes where refinement is made by subdividing an originally conforming mesh, the conforming mesh satisfies \eqref{eq:GeneralCondition} and hence Condition (F). If the polynomials used to approximate the metrics on the parent element are used within its children, Condition (F) also holds on the children. Clearly, the inherited metrics will differ from those computed directly from the element. This situation is natural in adaptive mesh refinement (AMR) where a conforming mesh is subdivided locally to enhance solution accuracy.
\item Half Order Elements\\
If $\spacevec X\in\mathbb{P}^{N/2}$ then the product ${{X_l}{\nabla _\xi }{X_m}}\in \mathbb{P}^{N}$ and the interpolations are exact on both the child and large elements. Thus, if the geometry is approximated at half the order of the solution (or less) everywhere, then the DGSEM is free-stream preserving. This approach is flexible in that one can generate a grid with a geometry half the order of the target solution order. Alternatively, given the geometry and mesh, the solution order can be chosen to be twice the order of the geometry.
\item Overintegration \\
Freestream preservation on a watertight non-conforming mesh with an mapping degree $N_g>N/2$ would be achieved if the numerical integration of the volume and surface terms are over-integrated \cite{Mengaldo201556},\cite{Kopriva2018} with a quadrature rule that exactly integrates both the basis of degree $N$ and the metrics of degree $2N_g$. 
\item \modColor{Transfer ${{\Gamma_l}{\nabla _\xi }{\Gamma_m}}$}\\
\modColor{In a recent approach presented in \cite{Kozdon2018}, the condition \eqref{eq:GeneralCondition} is enforced by a specific projection of the product ${{\Gamma_l}{\nabla _\xi }{\Gamma_m}}$ onto a finite element space of mixed degree $N-1$ and $N$ that is made continuous across non-conforming edges and surfaces. The metric terms are then computed via the discrete curl \eqref{eq:PolynomialMetrics}. If Gauss-Lobatto points are used for the approximation, then the surface differentiations will agree, modulo the scaling factor. Condition (V) will still hold via the curl relation. }
\end{enumerate}

\begin{rem}
The most general approaches for approximating the metric terms for non-conforming elements in three space dimension are numbers \modColor{six through eight}. The others are special cases.
\end{rem}
\section{Examples}
In this section, we present numerical evidence for free-stream preservation techniques cataloged in Sec. \ref{sec:catalog}. We simulate the Euler equations in three-dimensional fully periodic domains, and look at the errors against the free-stream solution. As shown in Fig.~\ref{fig:NonConformingPressFluc}, we encounter pressure fluctuations if the same polynomial degree for the solution ansatz and the element mappings are used. 

\begin{figure}[!htbp] 
   \centering
   \includegraphics[width=0.45\textwidth]{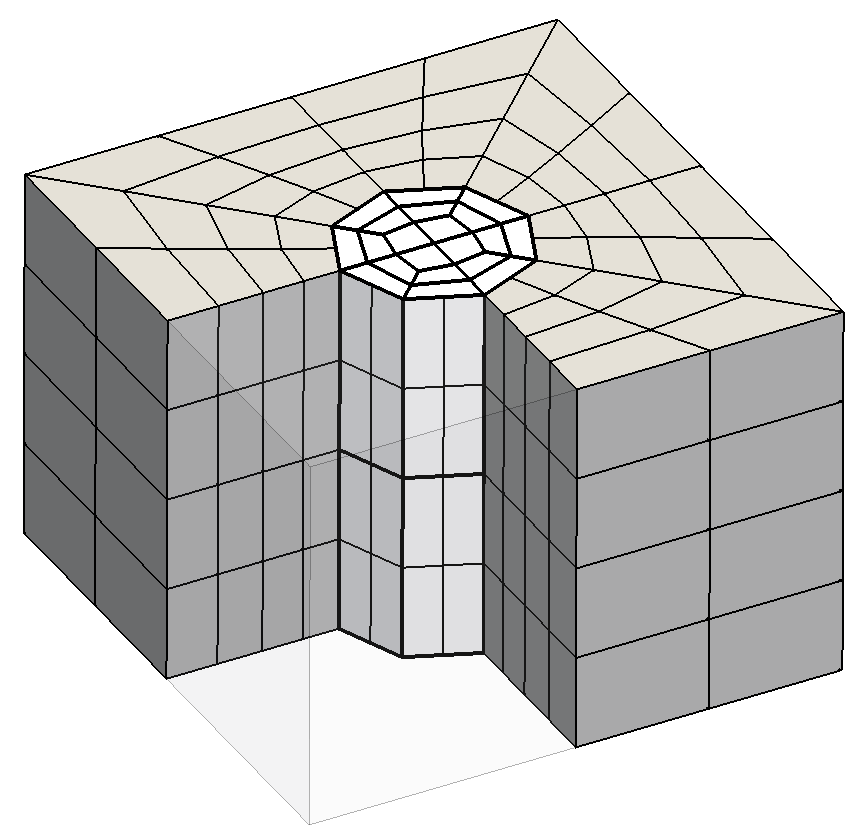} 
   \includegraphics[width=0.45\textwidth]{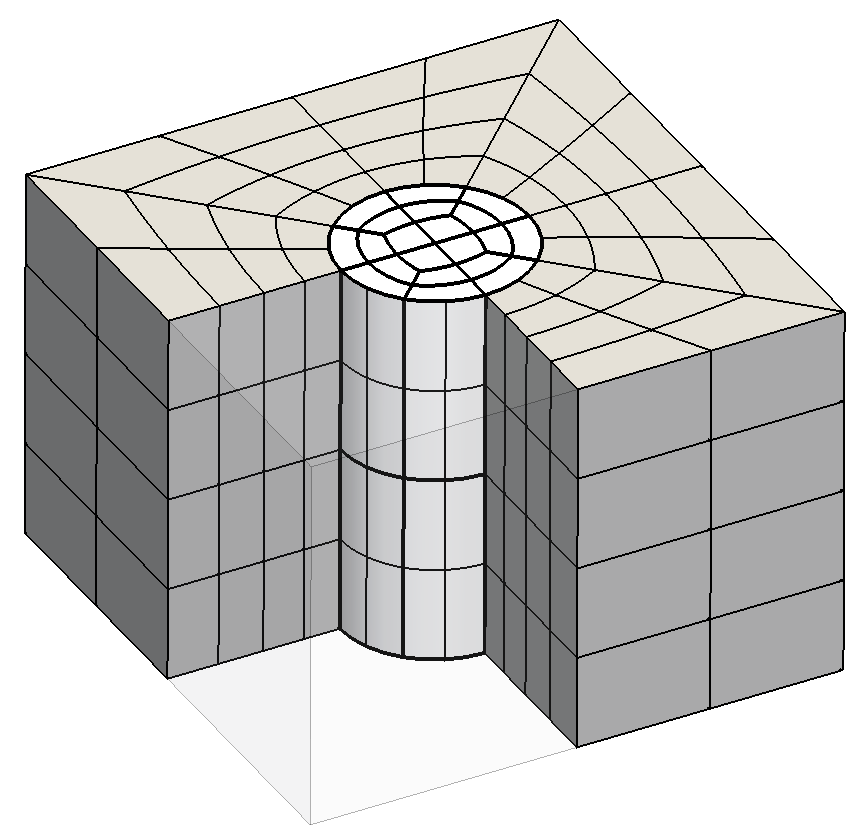} 
   \caption{`Cylinder' mesh: Extruded from two-dimensional linear ($N_g\!=\!1$) and curved mesh ($N_g\!=\!4$), the domains inside and outside the cylinder are conforming (different greyscale) and are connected with non-conforming interfaces.}
   \label{fig:cylindermeshes}
\end{figure}

\begin{figure}[!htbp] 
   \centering
   \includegraphics[width=0.45\textwidth]{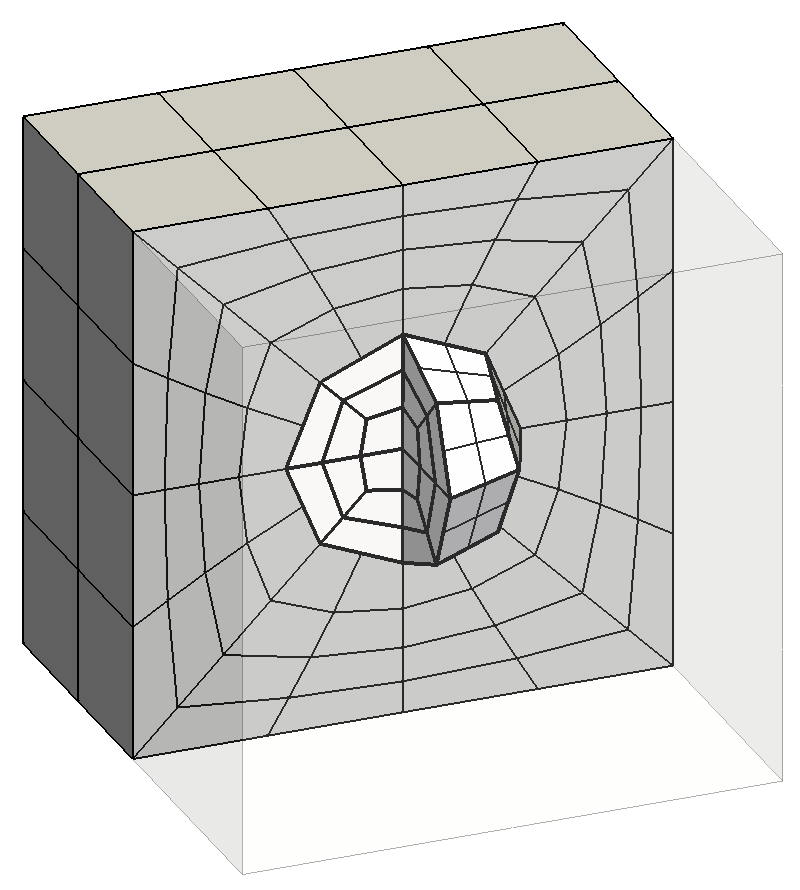} 
   \includegraphics[width=0.45\textwidth]{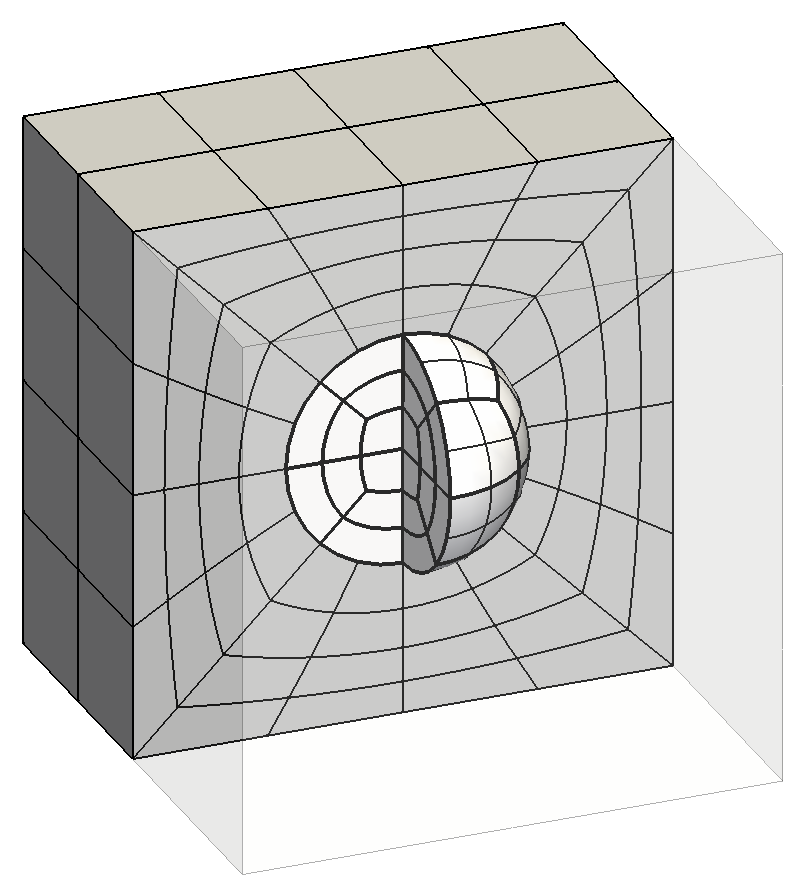} 
   \caption{`Sphere' mesh: Linear ($N_g\!=\!1$) and curved mesh ($N_g\!=\!4$), the domains inside and outside the sphere are conforming (different greyscale) and are connected with non-conforming interfaces.}
   \label{fig:spheremeshes}
\end{figure}

We investigate two choices of non-conforming meshes. One is extruded from a two-dimensional mesh, shown in Fig.~\ref{fig:cylindermeshes}, which we will call the `cylinder' mesh, since the non-conforming interfaces form a cylinder shape. The non-conforming interfaces in the second mesh form a sphere shape, thus called the `sphere' mesh which is three-dimensionally unstructured, as shown in Fig.~\ref{fig:spheremeshes}.  


\begin{figure}[!htbp] 
   \centering
   \includegraphics[width=0.49\textwidth]{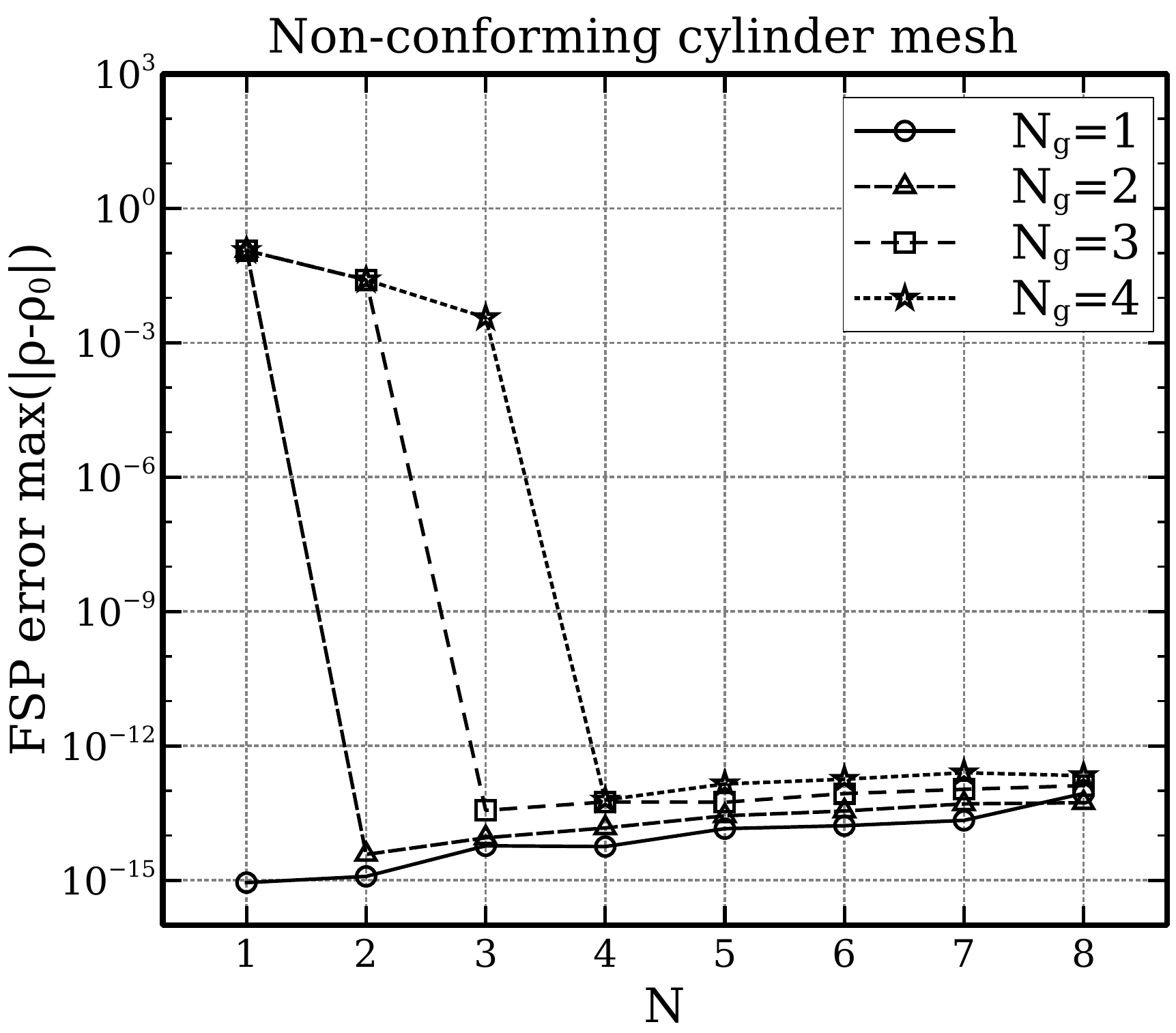} 
   \includegraphics[width=0.49\textwidth]{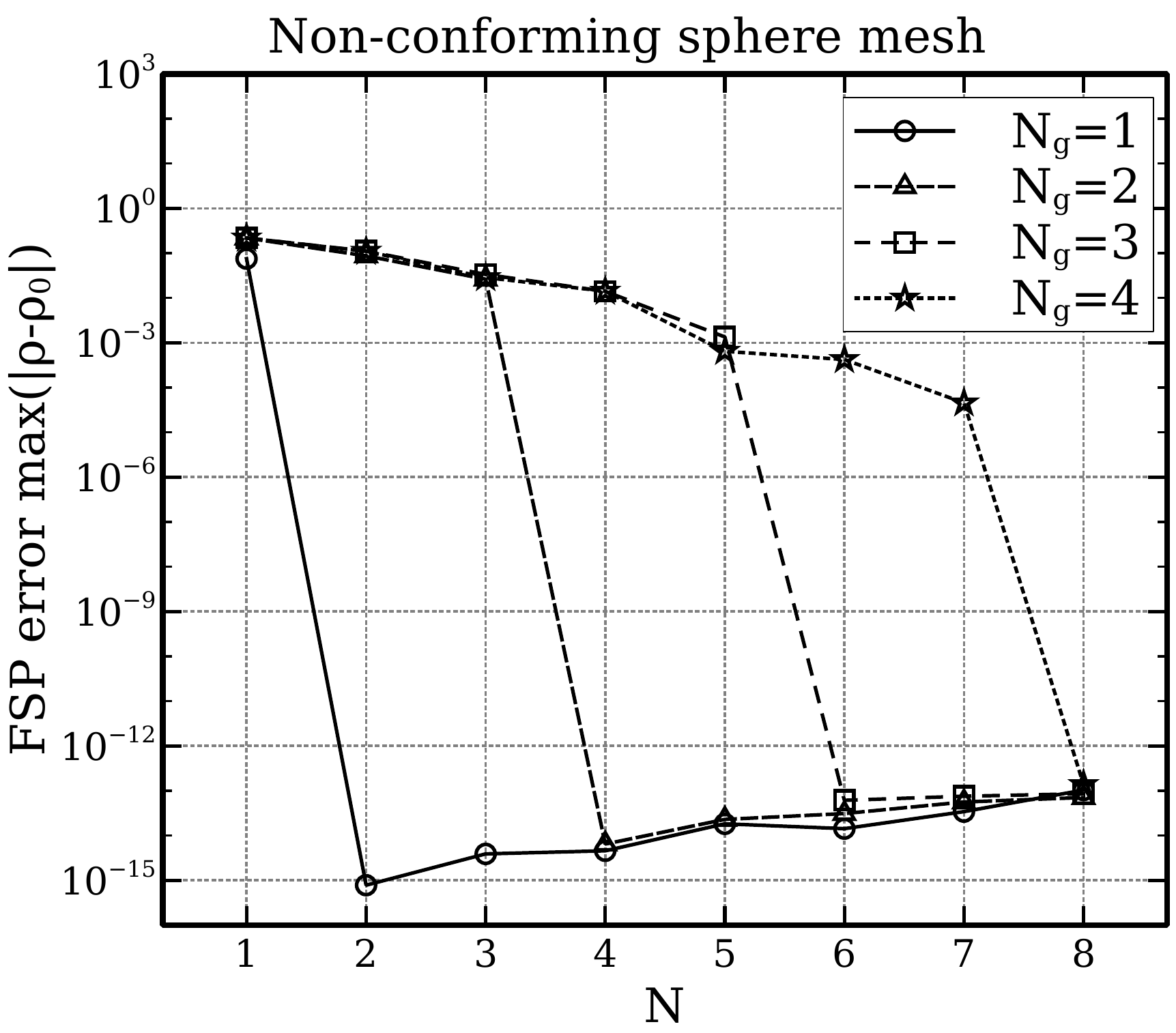} 
   \caption{Error of freestream preservation (density), for non-conforming cylinder extruded mesh and sphere mesh, for different degrees of the element mapping ($N_g$) over the polynomial degree of the solution ($N$).}
   \label{fig:FSPresults}
\end{figure}

All meshes were generated with the open-source tool HOPR\footnote{\url{https://github.com/fhindenlang/hopr} , \url{https://www.hopr-project.org} } (high order preprocessor) \cite{HOPR}. In HOPR, the mesh is built from straight-sided block definitions, with high order Gauss-Lobatto nodes in each element. Each interpolation point is then mapped via a smooth global mapping function to a curved geometry. Blocks can have non-conforming 2--1 or 4--1 interfaces. Here the `watertight condition' is enforced, which guarantees that at a non-conforming interface, the child face mappings are interpolations of the parent face mapping.  

We present the results of the freestream preservation tests on the non-conforming extruded cylinder mesh and the sphere mesh in Fig.~\ref{fig:FSPresults}. The simulations are run with Gauss nodes and the Lax-Friedrichs numerical surface flux. We initialize with the constant freestream $(\rho,v_1,v_2,v_3,p)=(0.7,0.2,0.3,-0.4,1.0)$ and run with periodic boundary conditions until a final simulated time of $T=1.0$. 

In Fig.~\ref{fig:FSPresults} we plot the maximum error norm of the density as a measure for freestream preservation, for increasing polynomial degree of the element mapping $N_g=1,\dots4$ over the polynomial degree of the solution $N=1,\dots,8$. The results demonstrate the observations of Sec. \ref{sec:catalog} that the freestream is preserved to machine rounding error on non-conforming meshes under the following conditions: For three-dimensional but extruded meshes, the polynomial degree of the solution must be at least equal to the degree of the element mapping ($N\geq N_g$), whereas in general three-dimensional geometry, we need the solution approximation order must be at least \emph{twice} the degree of the element mapping ($N\geq 2N_g$).
Note that if the conditions are not met, the freestream error is not negligible, irrespective of its spectral convergence seen in Fig.~\ref{fig:FSPresults}. It was also observed that these errors are not advected away but rather are constantly produced at the non-conforming interfaces. 

We found that the same conditions hold if we run the simulations with Gauss-Lobatto nodes instead of Gauss nodes. We also found that the scheme with Gauss-Lobatto nodes was unstable for $1<N<N_g$, though this choice is not relevant in practice since the watertight condition on the mesh cannot be met. 


\section{Summary}

The under-integration of the surface terms \cite{Mengaldo201556},\cite{Kopriva2018}, in the discontinuous Galerkin approximation introduces errors at the element faces that do not cancel and lead to free-stream preservation errors. These errors are not necessarily small, though the spectral polynomial approximation will at least ensure that the error decays exponentially fast. The general condition for free-stream preservation requires both volume (Condition V, eq. \eqref{eq:VolumeMetricIdentity}) and face (Condition F, eq. \eqref{eq:SurfacemetricIdentityII}) conditions. Condition V must hold for both conforming \cite{kopriva2006} and non-conforming meshes. Condition F, which is easily satisfied for conforming meshes, requires that the difference between the metric terms over the faces be orthogonal to the approximation space. We have cataloged eight sufficient conditions on the geometry that do not require the solution of a least squares (projection) solution to ensure Condition F. The most general are to approximate the geometry to no more than half the order of the solution, overintegration of quadratures, or by transferring the product ${{\Gamma_l}{\nabla _\xi }{\Gamma_m}}$ from the parent to the child faces.

\acknowledgement{ This work was supported by a grant from the Simons Foundation (\#426393, David Kopriva). G.G and T.B. have been supported by the European Research Council (ERC) under the European Union's Eights Framework Program Horizon 2020 with the research project \textit{Extreme}, ERC grant agreement no. 714487. DAK would like to thank Mr. Andres Rueda for his helpful comments during the preparation of this paper.

\begin{appendix}
\section{Appendix}
For ensuring the watertight condition, we require the geometry interpolants of the large face and the child faces, shown in Fig.~\ref{fig:MultidomainInterpolation}, to match.
The interpolation of the interpolant from the large face to the child faces is exact. The interpolation of the product of two interpolants is not. For completeness, we present a proof here, though only for edges in a two dimensional mesh. The proof for two dimensional faces in a three dimensional mesh will follow the same approach due to the tensor product approximation.

\begin{figure}[tbp] 
   \centering
   \includegraphics[width=3.5in]{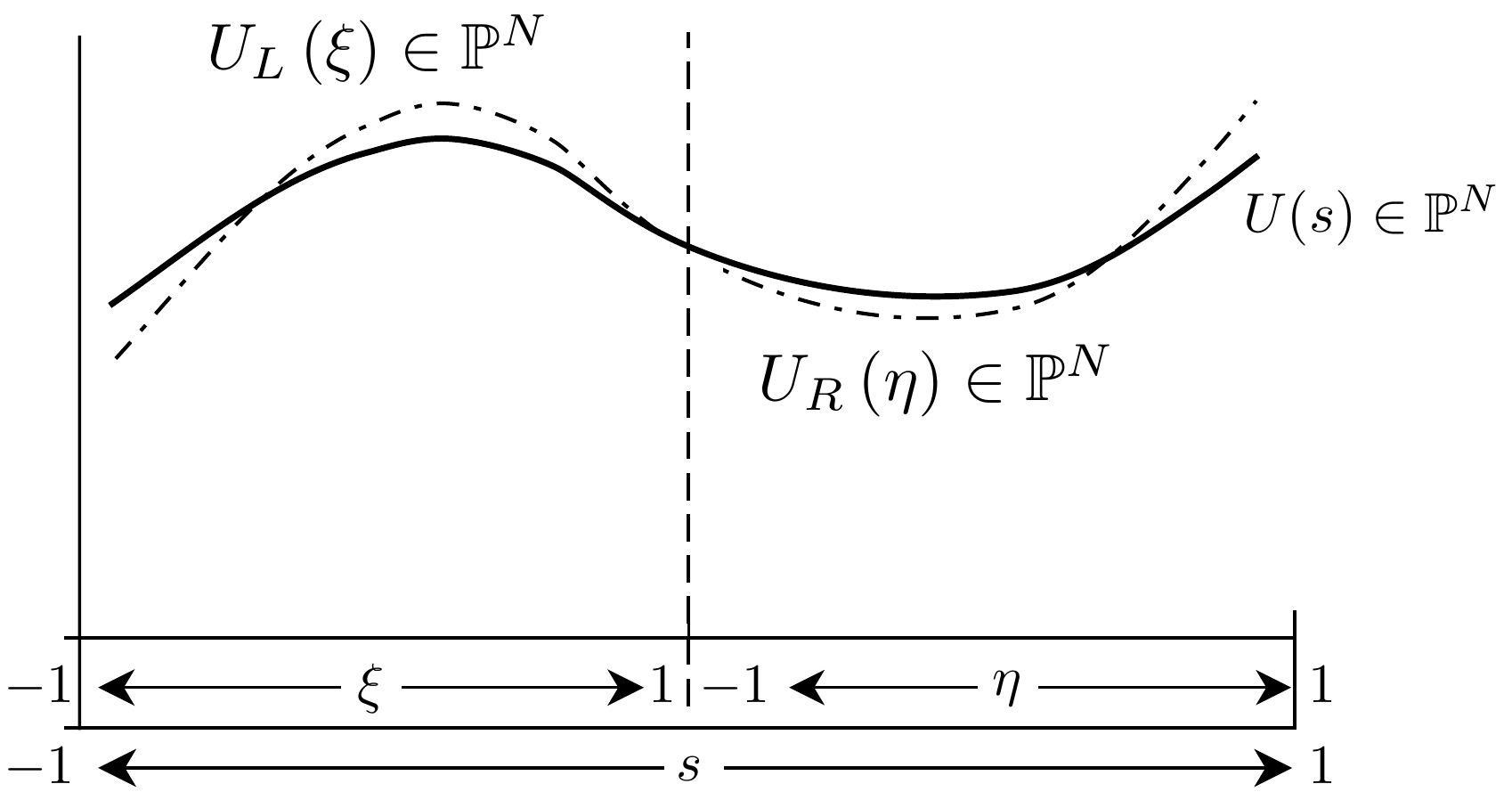} 
   \caption{Subdivision of interpolants. The polynomials $U^{L}$ and $U^{R}$ will be constructed to match $U$ on their respective intervals.}
   \label{fig:MultidomainInterpolation}
\end{figure}
Define the polynomial of degree $N$ over the interval $[-1,1]$ to be
\begin{equation}U\left( s \right) = \sum\limits_{j = 0}^N {{U_j}{\ell _j}\left( s \right)}. \quad s \in [-1,1]\end{equation}
Next subdivide the interval into two pieces $\xi \in [-1,1]$ and $\eta\in [-1,1]$ and two mappings
\begin{equation}s = \left\{ \begin{gathered}
  \frac{\xi-1 }{2}\quad \xi  \in [-1,1] \hfill \\
   \frac{\eta +1}{2}\quad \eta  \in [-1,1]. \hfill \\ 
\end{gathered}  \right.\end{equation}
Define the two polynomials that take on the values of $U$ at the interpolation points on the left and right intervals, i.e
\begin{equation}
  {U_{L}}\left( \xi  \right) = \sum\limits_{n = 0}^N {U\left( {\frac{{{\xi _n-1}}}{2}} \right){\ell _n}\left( \xi  \right)}  , \quad
  {U_{R}}\left( \eta  \right) = \sum\limits_{n = 0}^N {U\left( {\frac{{{\eta _n+1}}}{2}} \right){\ell _n}\left( \eta  \right)}.
  \label{eq:LRSplit}
 \end{equation}
Now we assume that we have two polynomials that match on the short and the long intervals. That is,
\begin{equation}
  U\left( s \right) = \sum\limits_{j = 0}^N {{U_j}{\ell _j}\left( s \right)} \quad U \ne 0, \quad
  V\left( s \right) = \sum\limits_{j = 0}^N {{V_j}{\ell _j}\left( s \right)} \quad V \ne 0.
\end{equation}
We then project the product of the two onto the polynomial space
\begin{equation}W(s) = {\mathbb{I}^N}\left( {UV} \right) = \sum\limits_{j = 0}^N {{U_j}{V_j}{\ell _j}\left( s \right)}, \end{equation}
and break the interval into two, as in \eqref{eq:LRSplit}, and define interpolants on each half
\begin{equation}\begin{gathered}
  {W_L}\left( \xi  \right) = \mathbb{I}_L^N\left( {UV} \right) = \sum\limits_{j = 0}^N {U\left( {\frac{{{\xi _j} - 1}}{2}} \right)V\left( {\frac{{{\xi _j} - 1}}{2}} \right){\ell _j}\left( \xi  \right)}  \hfill \\
  {W_R}\left( \xi  \right) = \mathbb{I}_R^N\left( {UV} \right) = \sum\limits_{j = 0}^N {U\left( {\frac{{{\eta _j} + 1}}{2}} \right)V\left( {\frac{{{h_j} + 1}}{2}} \right){\ell _j}\left( \eta  \right)} . \hfill \\ 
\end{gathered} \end{equation}
\\
Then we prove that
\begin{equation}
  W\left( {\frac{{\xi  - 1}}{2}} \right) \ne {W_L}\left( \xi  \right),\quad
  W\left( {\frac{{\eta  + 1}}{2}} \right) \ne {W_R}\left( \eta  \right),
\end{equation}
that is, the interpolation of the product onto the polynomials of degree $N$ on each half does not equal the interpolant of the product over the whole interval.

To prove this, we examine the error of the interpolants. We know from basic numerical analysis that
\begin{equation}E(s) = UV - W(s) = \frac{1}{{\left( {N + 1} \right)!}}\frac{{{\partial ^{N + 1}}\left( {UV} \right)}}{{\partial {s^{N + 1}}}}\prod\limits_{i = 0}^N {\left( {s - {s_i}} \right)} \end{equation}
and that $E(s)=0$ at precisely the $N+1$ interpolation points, $s_{i},\; i=0,1,\ldots,N$. Similarly for the interpolant on the left,
\begin{equation}{E_L}(\xi ) = U\left( {\frac{{\xi  - 1}}{2}} \right)V\left( {\frac{{\xi  - 1}}{2}} \right) - {W_L}(\xi ) = \frac{1}{{\left( {N + 1} \right)!}}\frac{{{\partial ^{N + 1}}\left( {UV} \right)}}{{\partial {\xi ^{N + 1}}}}\prod\limits_{i = 0}^N {\left( {\xi  - {\xi _i}} \right)}. \end{equation} That interpolant vanishes at precisely $N+1$ points $\xi_{i}$ and nowhere else unless the derivative of the product vanishes. 

The Gauss and Gauss-Lobatto points are symmetric and distinct about the middle of the interval and so there are twice as many (distinct) points $\xi_{i}$ on the interval $s\in[-1,0]$ as there are nodes $s_{i}$. Therefore, there exist nodes $\xi_{i}$ that are not equal to any node $s_{i}$. Let us choose, then, one such $\xi_{i}$ such that $s=(\xi_{i}-1)/2$ is not a node of the interpolation on $s$. Then $E_{L}\left(\xi_{i}\right)=0$, i.e., $U\left((\xi_{i}-1)/2\right)V\left((\xi_{i}-1)/2\right)=W_{L}\left(\xi_{i}\right)$, but $E\left((\xi_{i}-1)/2\right)\ne 0$ so that $U\left((\xi_{i}-1)/2\right)V\left((\xi_{i}-1)/2\right)\ne W\left((\xi_{i}-1)/2\right)$ and the result follows.

\end{appendix}

\bibliographystyle{plain}

\bibliography{dakBib}

\begin{thebibliography}{10}

\bibitem{Allaneau:2011pz}
Y.~Allaneau and Antony Jameson.
\newblock Connections between the filtered discontinuous {G}alerkin method and
  the flux reconstruction approach to high order discretizations.
\newblock {\em Computer Methods in Applied Mechanics and Engineering},
  200(49-52):3628--3636, 2011.

\bibitem{ISI:000306588600006}
Tormod Bjontegaard, Einar~M. Ronquist, and Oystein Trasdahl.
\newblock {Spectral Approximation of Partial Differential Equations in Highly
  Distorted Domains}.
\newblock {\em {Journal Of Scientific Computing}}, {52}({3}):{603--618}, {SEP}
  {2012}.

\bibitem{Bui-Tanh:2012zp}
Tan Bui-Tanh and Omar Ghattas.
\newblock Analysis of an $hp$-nonconforming discontinuous {G}alerkin spectral
  element method for wave propagation.
\newblock {\em SIAM J. Numerical Analysis}, 50(3):1801--1826, 2012.

\bibitem{canuto2007}
C.~Canuto, M.~Hussaini, A.~Quarteroni, and T.~Zang.
\newblock {\em Spectral Methods: Evolution to Complex Geometries and
  Applications to Fluid Dynamics}.
\newblock Springer, Berlin, 2007.

\bibitem{Fortunato2015}
Meire Fortunato and Per-Olof Persson.
\newblock High-order unstructured curved mesh generation using the winslow
  equations.
\newblock {\em Journal of Computational Physics},307:1-14, 2016.

\bibitem{Gassner2018}
Gregor~J. Gassner, Andrew~R. Winters, Florian~J. Hindenlang, and David~A.
  Kopriva.
\newblock The br1 scheme is stable for the compressible {N}avier--{S}tokes
  equations.
\newblock {\em Journal of Scientific Computing}, 77(1):154-200,2018.

\bibitem{Gassner:2016ye}
Gregor~J Gassner, Andrew~R Winters, and David~A Kopriva.
\newblock Split form nodal discontinuous {G}alerkin schemes with
  summation-by-parts property for the compressible {E}uler equations.
\newblock {\em {Journal Of Computational Physics}}, 327:39--66, 2016.

\bibitem{Gordon&Hall1973a}
W.J. Gordon and C.A. Hall.
\newblock Construction of curvilinear co-ordinate systems and their
  applications to mesh generation.
\newblock {\em International Journal for Numerical Methods in Engineering
  Engineering}, 7:461--477, 1973.

\bibitem{Hestahven:1008th}
J.S. Hesthaven and T.~Warburton.
\newblock {\em Nodal Discontinuous {G}alerkin Methods: Algorithms, Analysis,
  and Applications}.
\newblock Springer, 2008.

\bibitem{HOPR}
F.~Hindenlang, T.~Bolemann, and C-D. Munz.
\newblock Mesh curving techniques for high order discontinuous {G}alerkin
  simulations.
\newblock In {\em IDIHOM: Industrialization of High-Order Methods-A Top-Down
  Approach}, pages 133--152. Springer, 2015.

\bibitem{Hindenlang:2014gl}
Florian Hindenlang.
\newblock {\em Mesh Curving Techniques for High Order Parallel Simulations on
  Unstructured Meshes}.
\newblock PhD thesis, University of Stuttgart, 2014.

\bibitem{Johnen:2014eu}
A.~Johnen, J.-F. Remacle, and C.~Geuzaine.
\newblock Geometrical validity of high-order triangular finite elements.
\newblock {\em Engineering with Computers}, 30(3):375--382, 2014.

\bibitem{Karniadakis:2005fj}
George~Em Karniadakis and Spencer~J. Sherwin.
\newblock {\em Spectral/hp Element Methods for Computational Fluid Dynamics}.
\newblock Oxford University Press, 2005.

\bibitem{Kopera201492}
Michal~A. Kopera and Francis~X. Giraldo.
\newblock Analysis of adaptive mesh refinement for {IMEX} discontinuous
  {G}alerkin solutions of the compressible euler equations with application to
  atmospheric simulations.
\newblock {\em Journal of Computational Physics}, 275:92 -- 117, 2014.

\bibitem{Kopriva:1996:JCP96b}
D.~A. Kopriva.
\newblock A conservative staggered-grid {C}hebyshev multidomain method for
  compressible flows {II}. {A} semi-structured method.
\newblock {\em Journal of Computational Physics}, 128(2):475--488, 1996.

\bibitem{kopriva2006}
D.~A. Kopriva.
\newblock Metric identities and the discontinuous spectral element method on
  curvilinear meshes.
\newblock {\em Journal of Scientific Computing}, 26(3):301--327, 2006.

\bibitem{gassner2010}
D.~A. Kopriva and G.~Gassner.
\newblock On the quadrature and weak form choices in collocation type
  discontinuous {G}alerkin spectral element methods.
\newblock {\em Journal of Scientific Computing}, 44(2):136--155, 2010.

\bibitem{Koprivaetal2000}
D.~A. Kopriva, S.~L. Woodruff, and M.~Y. Hussaini.
\newblock Computation of electromagnetic scattering with a non-conforming
  discontinuous spectral element method.
\newblock {\em International Journal for Numerical Methods in Engineering},
  53(1):105--122, 2002.

\bibitem{Kopriva:2009nx}
David~A. Kopriva.
\newblock {\em Implementing Spectral Methods for Partial Differential
  Equations}.
\newblock Scientific Computation. Springer, May 2009.

\bibitem{Kopriva:2017yg}
David~A. Kopriva.
\newblock A polynomial spectral calculus for analysis of {DG} spectral element
  methods.
\newblock In Hesthaven~J. Bittencourt~M., Dumont~N., editor, {\em Spectral and
  High Order Methods for Partial Differential Equations ICOSAHOM 2016}, volume
  119, pages 21--40, Cham, 2017. Springer.

\bibitem{Kopriva2018}
David~A. Kopriva.
\newblock Stability of overintegration methods for nodal discontinuous
  {G}alerkin spectral element methods.
\newblock {\em Journal of Scientific Computing}, 76(1):426--442, Jul 2018.

\bibitem{Kopriva2016274}
David~A. Kopriva and Gregor~J. Gassner.
\newblock Geometry effects in nodal discontinuous {G}alerkin methods on curved
  elements that are provably stable.
\newblock {\em Applied Mathematics and Computation}, 272, Part 2:274 -- 290,
  2016.

\bibitem{ISI:A1986A176500008}
K.Z. Korczak and A.T. Patera.
\newblock An isoparametric spectral element method for solution of the
  {N}avier-{S}tokes equations in complex-geometry.
\newblock {\em {Journal Of Computational Physics}}, {62}({2}):{361--382}, {FEB}
  {1986}.

\bibitem{Kovalev:2005qv}
Konstantin Kovalev.
\newblock {\em Unstructured Hexahedral Non-conformal Mesh Generation}.
\newblock PhD thesis, Faculty of Engineering Vrije Universiteit Brussel
  Belgium, December 2005.

\bibitem{Kozdon2018}
Jeremy~E. Kozdon and Lucas~C. Wilcox.
\newblock An energy stable approach for discretizing hyperbolic equations with
  nonconforming discontinuous {G}alerkin methods.
\newblock {\em Journal of Scientific Computing}, 76(3):1742--1784, Sep 2018.

\bibitem{L.-Friedrichs:2017hm}
A.~R. Winters G. J. Gassner D. W.~Zingg L.~Friedrichs, D. C. D.
  R.~Fern{\'a}ndez and J.~Hicken.
\newblock Conservative and stable degree preserving sbp operators for
  non-conforming meshes.
\newblock {\em Journal of Scientific Computing},
  https://doi.org/10.1007/s10915-017-0563-z, 2017.

\bibitem{Lovgren:2009zr}
A.E. Lovgren and Y.~Maday an~E.M~Ronquist.
\newblock Global $c^1$ maps on general domains.
\newblock {\em Mathematical Models and Methods in Applied Sciences},
  19(5):803--832, 2009.

\bibitem{Lucas-Friedrich:2018ng}
David C. Del Rey Fernandez Gregor J. Gassner Matteo Parsani Mark H.~Carpenter
  Lucas~Friedrich, Andrew R.~Winters.
\newblock An entropy stable h/p non-conforming discontinuous {G}alerkin method
  with the summation-by-parts property.
\newblock {\em Journal of Scientific Computing}, arXiv:1712.10234, 2018.

\bibitem{Mengaldo2016}
G.~Mengaldo, D.~De~Grazia, P.~E. Vincent, and S.~J. Sherwin.
\newblock On the connections between discontinuous {G}alerkin and flux
  reconstruction schemes: Extension to curvilinear meshes.
\newblock {\em Journal of Scientific Computing}, 67(3):1272--1292, Jun 2016.

\bibitem{Mengaldo201556}
G.~Mengaldo, D.~De Grazia, D.~Moxey, P.E. Vincent, and S.J. Sherwin.
\newblock Dealiasing techniques for high-order spectral element methods on
  regular and irregular grids.
\newblock {\em Journal of Computational Physics}, 299:56 -- 81, 2015.

\bibitem{Moxey:2015xe}
D.~Moxey, M.~Hazan, S.J. Sherwin, and J.~Peiro.
\newblock Curvilinear mesh generation for boundary layer problems.
\newblock In Norbert Kroll, Charles Hirsch, Francesco Bassi, Craig Johnston,
  and Koen Hillewaert, editors, {\em IDIHOM: Industrialization of High-Order
  Methods - A Top-Down Approach}, volume 128 of {\em Notes on Numerical Fluid
  Mechanics and Multidisciplinary Design}, pages 41--64. Springer International
  Publishing, 2015.

\bibitem{NelsonDissertation}
Daniel A.~W. Nelson.
\newblock {\em High-Fidelity Lagrangian Coherent Structures Analysis and DNS
  with Discontinuous-{G}alerkin Methods}.
\newblock PhD thesis, University of California at San Diego, 2015.

\bibitem{Nonomura2010197}
Taku Nonomura, Nobuyuki Iizuka, and Kozo Fujii.
\newblock Freestream and vortex preservation properties of high-order {WENO}
  and {WCNS} on curvilinear grids.
\newblock {\em Computers \& Fluids}, 39(2):197 -- 214, 2010.

\bibitem{Patera:1984:JCP}
A.T. Patera.
\newblock A spectral element method for fluid dynamics - laminar flow in a
  channel expansion.
\newblock {\em Journal of Computational Physics}, 54(3):468--488, 1984.

\bibitem{Persson20091585}
P.-O. Persson, J.~Bonet, and J.~Peraire.
\newblock Discontinuous {G}alerkin solution of the {N}avier-{S}tokes equations
  on deformable domains.
\newblock {\em Computer Methods in Applied Mechanics and Engineering},
  198(17--20):1585 -- 1595, 2009.

\bibitem{Persson:2009mz}
Per-Olof Persson and Jaime Peraire.
\newblock Curved mesh generation and mesh refinement using lagrangian solid
  mechanics.
\newblock {\em Preprint}, 2009.

\bibitem{NME:NME2800}
Matthew~L. Staten, Jason~F. Shepherd, Franck Ledoux, and Kenji Shimada.
\newblock Hexahedral mesh matching: Converting non-conforming
  hexahedral-to-hexahedral interfaces into conforming interfaces.
\newblock {\em International Journal for Numerical Methods in Engineering},
  82(12):1475--1509, 2010.

\bibitem{Thomas&Lombard1979}
P.D. Thomas and C.K. Lombard.
\newblock Geometric conservation law and its application to flow computations
  on moving grids.
\newblock {\em AIAA Journal}, 17(10):1030--1037, 1979.

\bibitem{Visbal1999}
M.R. Visbal and D.V. Gaitonde.
\newblock High-order accurate methods for complex unsteady subsonic flows.
\newblock {\em AIAA Journal}, 37(10):1231--1239, 1999.

\bibitem{Visbal2002}
M.R. Visbal and D.V. Gaitonde.
\newblock On the use of higher-order finite-difference schemes on curvilinear
  and deforming meshes.
\newblock {\em Journal of Computational Physics}, 181:155--185, 2002.

\bibitem{Xie:2013rw}
ZhongQ. Xie, Ruben Sevilla, Oubay Hassan, and Kenneth Morgan.
\newblock The generation of arbitrary order curved meshes for 3d finite element
  analysis.
\newblock {\em Computational Mechanics}, 51(3):361--374, 2013.

\bibitem{Zhang2015147}
Bin Zhang and Chunlei Liang.
\newblock A simple, efficient, and high-order accurate curved sliding-mesh
  interface approach to spectral difference method on coupled rotating and
  stationary domains.
\newblock {\em Journal of Computational Physics}, 295:147 -- 160, 2015.

\end{thebibliography}

\end{document}